\documentstyle[amscd, amssymb, verbatim]{amsart}
\theoremstyle{plain}
\newtheorem{Thm}{Theorem}
\newtheorem{Cor}{Corollary}
\newtheorem{Prop}{Proposition}[section]
\newtheorem{Lem}[Prop]{Lemma}
\newtheorem{Def}{Definition}

\input{epsf}

\theoremstyle{definition}

\begin{document}
% this is the second version of the introduction!

%Sample file: article.tpl
%Typeset with AMSlaTex format file
%Preamble 
% Style section 

%\begin{document}
%Topmatter
\title
{On The Homflypt Skein Module of $ S^{1} \times S^2$}

\author{Patrick M. Gilmer}
\email{gilmer@@math.lsu.edu}
\author{Jianyuan Zhong}
\email{zhong@@math.lsu.edu}
\address{Department of Mathematics\\
   Louisiana State University\\
    Baton Rouge, LA 70803}

\date{June 28, 2000}

%End topmatter

\begin{abstract}
   Let $k$ be a subring of the field of rational functions in $x, v, s$ which contains $x^{\pm 1}, v^{\pm 1}, s^{\pm 1}$. If $M$ is an oriented $3$-manifold,  let $S(M)$ denote the Homflypt skein module of $M$ over $k$. This is the free $k$-module generated by isotopy classes of framed oriented links in $M$ quotiented by the Homflypt skein relations: (1) $x^{-1}L_{+}-xL_{-}=(s-s^{-1})L_{0}$; (2) $L$ with a positive twist $=(xv^{-1})L$; (3) $L\sqcup O=(\frac{v-v^{-1}}{s-s^{-1}})L$ where $O$ is the unknot. We give two bases for the relative Homflypt skein module of the solid torus with $2$ points in the boundary. The first basis is related to the basis of $S(S^1\times D^2)$ given by V. Turaev and also J. Hoste and M. Kidwell; the second basis is related to a Young idempotent basis for $S(S^1\times D^2)$ based on the work of A. Aiston, H. Morton and C. Blanchet. We prove that if the elements $s^{2n}-1$, for $n$ a nonzero integer, and the elements $s^{2m}-v^{2}$, for any integer $m$, are invertible in $k$, then $S(S^{1} \times S^2)=k$-torsion module $\oplus k$. Here the free part is generated by the empty link $\phi$. In addition, if the elements $s^{2m}-v^{4}$, for $m$ an integer, are invertible in $k$, then $S(S^{1} \times S^2)$ has no torsion. We also obtain some results for more general $k$.

\end{abstract}

\maketitle
\section{Introduction} 
The Kauffman bracket skein module of $S^1 \times S^2$ was discussed by J. Hoste and J. Przytycki in \cite{HP95a}. This motivated us to investigate the Homflypt skein module of $S^1 \times S^2$. We show this skein module, over $F[x, x^{-1}]$, where $F$ is the field of rational functions in $ v,\ s$, is free on one generator, the empty link.

Let $k$ be an integral domain containing the invertible elements $x$, $v$ and $s$. Moreover we assume that $s-s^{-1}$ is invertible in $k$. 

We will be working with framed oriented links. By this we mean links equipped with a string orientation together with a nonzero normal vector field up to homotopy. The links described by figures in this paper will be assigned the ``blackboard'' framing which points to the right when travelling along an oriented strand.

\begin{Def}
{\bf The Homflypt skein module}. Let $M$ be an oriented $3$-manifold. The Homflypt skein module of $M$, denoted by $S(M)$, is the $k$-module freely generated by isotopy classes of framed oriented links in $M$ including the empty link, quotiented by the Homflypt skein relations given in the following figure.
$$
x^{-1}\quad\raisebox{-3mm}{\epsfxsize.3in\epsffile{left.ai}}\quad
        -\quad x\quad\raisebox{-3mm}{\epsfxsize.3in\epsffile{right.ai}}\quad
=\quad (\ s - \ s^{-1})\quad\raisebox{-3mm}{\epsfxsize.3in\epsffile{parra.ai}}\quad ,
$$
$$
\raisebox{-3mm}{\epsfxsize.35in\epsffile{framel.ai}}
        \quad=\quad (xv^{-1})\quad \raisebox{-3mm}{\epsfysize.3in\epsffile{orline.ai}}\quad  ,
$$
$$
L\ \sqcup \raisebox{-2mm}{\epsfysize.3in\epsffile{unknot.ai}}\quad
                                        =\quad {\dfrac{v^{-1}-v} {\ s - \ s^{-1}}}\quad L\quad.
$$
\end{Def} 

An embedding $f: M\to N$ induces a well defined homomorphism $f_{*}: S(M)\to S(N)$. Since $S^1 \times S^2$ can be obtained by adding a $2$-handle and a $3$-handle to the solid torus $S^1\times D^2$, we have an induced surjective map of the skein modules $S(S^1 \times D^2)\to S(S^1 \times S^2)$. V.~Turaev \cite{vT} gave a linear basis for $S(S^1 \times D^2)$  whose elements are represented by the collection of the monomials $A_{i_1}A_{i_2}\cdots A_{i_j}$ in the commuting elements $A_1, A_{-1},  A_2, A_{-2}, \dotsc$. We include the empty monomial, which is denoted by $1$ and occasionally by $A_{0}$, and represents the empty diagram. If $j=0$, the expression $A_{i_1}A_{i_2}\cdots A_{i_j}$ denotes $1$. $A_1$ represents the closure of the $1$-string braid oriented in the clockwise direction. The element $A_j$ represents the closure of the $j$-string braid $\sigma_{j-1}\sigma_{j-2}\ldots\sigma_1$ with strings oriented in the clockwise direction; and $A_{-j}$ represents the closure of the $j$-string braid $\sigma_{j-1}^{-1}\sigma_{j-2}^{-1}\ldots\sigma_{1}^{-1}$ with strings oriented in the counterclockwise direction. Here we follow the Morton-Aiston convention for $\sigma_{i}$, where $\sigma_{i}$ is the positive permutation braid corresponding to the transposition $(i\ i+1)$. A diagram representing $A_{1}A_{2}A_{-3}$ is shown below. Independently J. Hoste and M. Kidwell \cite{HK90} found the same description for $S(S^1 \times D^2)$. We will refer to the above basis as the monomial basis  for $S(S^1 \times D^2)$.

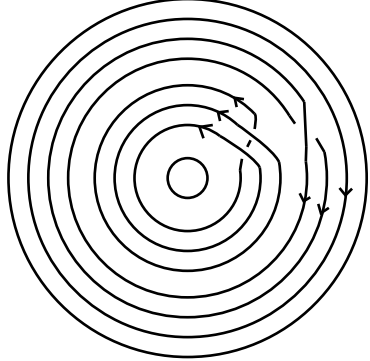
\hfil {\epsffile{a1a2a3.ai}} \hfil

$S(S^1 \times D^2)$ forms a commutative algebra with multiplication induced by embedding two solid tori in a single solid torus in a standard way. Let $C^{+}$ be the subalgebra freely generated by $\{A_{j} \mid j \in {\Bbb Z}, j\geq 0 \}$.  In other words, $C^{+}$ is the subspace of $S(S^1 \times D^2)$ generated by the closure of braids with the clockwise orientation and the empty link.  For $n>0$, let $C_{n}$ be the submodule generated by $\{A_{i_1}A_{i_2}\cdots A_{i_j} \mid \sum_{m=1}^{j} i_m=n, i_m \in {\Bbb Z}, i_m> 0, 1\leq m\leq j\}$, and $C_{0}$ is the submodule generated by $A_{0}$, which is the empty link. 
 
As a linear space $C^{+}$ is graded by
$$C^{+} \cong  \bigoplus_{n\geq 0} C_n .$$

Let $C^{-}$ be the subalgebra freely generated by $\{A_{-j} \mid j \in {\Bbb Z}, j\geq 0 \}$.  In other words, $C^{-}$ is the subspace of $S(S^1 \times D^2)$ generated by the closure of braids with the counterclockwise orientation and the empty link. For $n>0$, let $C_{-n}$ be the submodule generated by $\{A_{-i_1}A_{-i_2}\cdots A_{-i_j} \mid \sum_{m=1}^{j} i_m=n, i_m \in {\Bbb Z}, i_m> 0, 1\leq m\leq j\}$. 
 
As a linear space $C^{-}$ is graded by
$$C^{-} \cong  \bigoplus_{m\geq 0} C_{-m} .$$
Moreover, we have a module decomposition as
\begin{equation}
S(S^1\times D^2)\cong\bigoplus_{m \geq 0, n \geq 0}C_{-m} \otimes C_{n}  \tag{1.1}
\end{equation}

If $s^{2n}-1$ is invertible for integers $n > 0$, A. Aiston and H. Morton presented a new basis for $C^{+}$. The basis elements are indexed by the Young diagrams. For each Young diagram $\lambda$, one has a basis element $Q_{\lambda}$. See diagram.
$$Q_{\lambda}=\raisebox{-20mm}{\ \epsfxsize0in
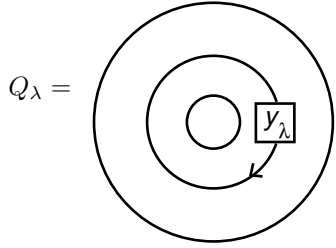
\epsffile{ylaa.ai}}$$
Here we draw one single string in the picture for $|\lambda|$ parallel strings in the shape of $\lambda$, where $|\lambda|$ is the size of the Young diagram $\lambda$, and the box labelled by $y_{\lambda}$ represents a certain linear combination of braid diagrams associated to $\lambda$. We call $y_{\lambda}$ the Young idempotent corresponding to the Young diagram $\lambda$, the definition of $y_{\lambda}$ is in section 3. The following proposition follows from the Morton-Aiston result above and Equation (1.1).

\begin{Prop}
If $s^{2n}-1$ is invertible for integers $n > 0$, $S(S^1 \times D^2)$ has a countable infinite basis given by $Q_{\lambda, \mu}$, where $\lambda, \ \mu$ vary over all Young diagrams. See diagram. 
\end{Prop}
$$Q_{\lambda, \mu}=\raisebox{-20mm}{\ \epsfxsize0in
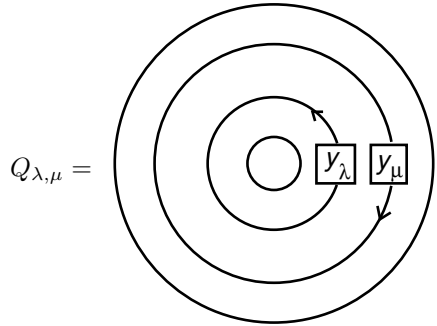
\epsffile{ylamu.ai}}$$
%Here $Q_{0,0}$ denotes the empty link, we will denote it by $Q_{0, \mu}$ if %$|\lambda|=0$, and $Q_{1, \mu}$ if $|\lambda|=1$; both $y_{\lambda}$ and %$y_{\mu}$ are linear combinations of the monomial basis elements. 

\noindent We will refer to this basis as the Young idempotent basis.

Now we consider the space $S^1 \times S^2$ as obtained by adding a 2-handle and a 3-handle to the solid torus. The addition of the $2$-handle will result in relations between the generators $Q_{\lambda, \mu}$. G.~Masbaum \cite{gM96} makes use of the relative Kauffman bracket skein of $S^1 \times D^2$ with two points in the boundary to recover and refine J. Hoste and J. Przytycki's calculation of the Kauffman bracket skein modules of $S^1 \times S^2$ and lens spaces; for the same reason,  we study the relative Homflypt skein module of $S^1 \times D^2$ with two points in the boundary. We use it to parametrize the relations. If $M$ has nonempty boundary, a framed point $y$ in $\partial M$ is a point together with a vector based at $y$ which is tangent to $\partial M$.

\begin{Def}
{\bf The relative Homflypt skein module}. Let $X=\{x_1, x_2,\cdots, x_{n}\}$ be a finite set of framed points oriented negatively (called input points) in $\partial M$, and let $Y=\{y_1, y_2,\cdots, y_{n}\}$ be a finite set of framed points oriented positively (called output points) in the boundary $\partial M$. Define the relative skein module $S(M, X, Y)$ to be the $k$-module generated by relative framed oriented links in $(M,\partial M)$ such that $L \cap \partial M=\partial L=\{x_i, y_i\}$ with the induced framing and orientation, considered up to an ambient isotopy fixing $\partial M$, quotiented by the Homflypt skein relations. 
\end{Def}

In particular, we will study the relative Homflypt skein module of the solid torus with an input point $A$ and an output point $B$. By a slight abuse of notation, we denote this module by $S(S^1 \times D^2, A, B)$. 

We will give two bases for $S(S^1 \times D^2, A, B)$. One is related to the monomial basis of $S(S^1 \times D^2)$ and is described in section $2$. The other basis for $S(S^1 \times D^2, A, B)$ is related to the Young idempotent basis and is given by the following theorem. 

\begin{Thm}
If $s^{2n}-1$ is invertible for integers $n> 0$, $S(S^1 \times D^2, A, B)$ has a countable infinite basis given by the collection of elements of $Q_{\lambda, \mu,c}', Q_{\lambda, \mu, c}'' $, where $\lambda$ and $\mu$ vary over all Young diagrams and $c$ varies over all extreme cells of $\mu$.  
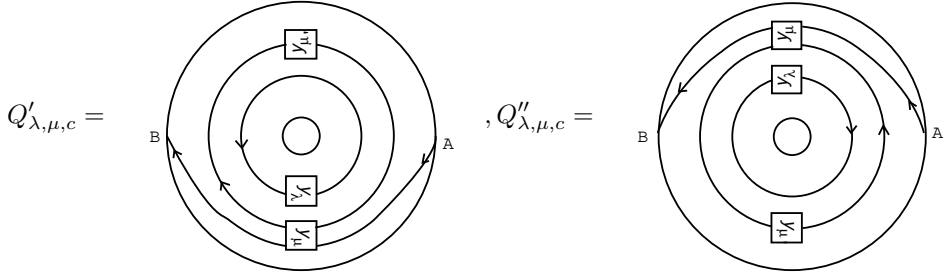
$$
Q_{\lambda, \mu,c}'=\raisebox{-20mm}{\ \epsfxsize1.9in\epsffile{abr1.ai}},
Q_{\lambda, \mu,c}''=\raisebox{-20mm}{\ \epsfxsize1.9in\epsffile{abr2.ai}}
$$
\end{Thm}
%$$
%Q_{\lambda, \mu,c}'=\raisebox{-20mm}{\ \epsfxsize1.6in\epsffile{ylamu1a.ai}},
%Q_{\lambda, \mu,c}''=\raisebox{-20mm}{\ \epsfxsize1.6in\epsffile{ylamu3.ai}}
%$$

\noindent Here the Young diagram $\mu '$ is obtained from the Young diagram $\mu$ by removing the extreme cell $c$, where an extreme cell is a cell such that if we remove it, we obtain a legitimate Young diagram. We prove this theorem in section 4. 

Let $\phi$ denote the empty link in $S^{1} \times S^2$, and $<\phi>$ denote the $k$-submodule of $ S(S^{1} \times S^2)$ generated by $\phi$.

\begin{Thm} If $s^{2n}-1$ is invertible for $n> 0$, $n\in {\Bbb Z}$, and $m\in {\Bbb Z}$,

(i) $ S(S^{1} \times S^2)\ /<\phi>$ is a $k$-torsion module. 

(ii) Assuming all elements of the form $s^{2m}-v^{2}$ are invertible in $k$, then $S(S^{1} \times S^2)={\hbox{k-torsion module}}\  \oplus <\phi>$.

(iii) Assuming in addition that all elements of the form $s^{2m}- v^{4}$ are invertible in $k$, then $S(S^{1} \times S^2)=<\phi>$.
\end{Thm}
\noindent We do not know whether these torsion submodules are nonzero. We prove this theorem in section 5 and the following theorem in section 6. 
%\begin{Cor}
%Let $k={\cal F}[x, x^{-1}]$, where ${\cal F}$ is a field of rational functions %in $v, s$, then $S(S^{1} \times S^2)=<\phi>$.
%\end{Cor}

\begin{Thm}
Let $k$ be a subring of the field of rational functions in $x, v, s$ which contains $x^{\pm 1}, v^{\pm 1}, s^{\pm 1}$, then $<\phi>$ is a free submodule of $ S(S^{1} \times S^2)$.
\end{Thm}

Let ${\cal {R}}$ be the subring of the field of rational functions in $v$, $s$ generated by $v$, $s$, $v^{-1}$, $s^{-1}$, $\dfrac{1}{s^{2n}-1}$ for $n>0$ and $\dfrac{1}{s^{2m}-v^4}$ for all $m$. 

\begin{Cor}
Let $k={\cal R}[x, x^{-1}]$, then $S(S^{1} \times S^2)=<\phi>$ is a free module generated by $\phi$.
\end{Cor}

Our work is clarified by the following proposition:
\begin{Prop}
Let $M$ be an oriented $3$-manifold, and let $H_{1}(M)$ denote the first homology group of $M$, then 
$$S(M)=\bigoplus_{z \in H_{1}(M)} S_{z}(M).$$
\end{Prop}
\noindent Here $S_{z}(M)$ is the submodule generated by the isotopy classes of framed oriented links in $M$ representing the homology class $z$, quotiented by the Homflypt skein relations.

\begin{pf}
This follows from the fact that the Homflypt skein relations respect homology classes.
\end{pf}

Thus, the modules $S(S^1\times D^2)$ and $S(S^1\times S^2)$ are ${\Bbb Z}$-graded. It is interesting to note that as a graded and commutative algebra, $S(S^1\times D^2) \cong C^{+} \otimes C^{-}$.

{\em Remark.} There is also a relative version of Proposition $1.2$.

\begin{Prop}
If $s^{2n}-1$ is invertible for $n> 0$, and $s^{2m}-v^{2}$ is invertible for all $m$, then $S_{0}(S^{1} \times S^2)=<\phi>$.
\end{Prop}
\noindent Here $S_{0}(S^{1} \times S^2)$ is the submodule generated by the isotopy classes of framed oriented null homologous links.

Let ${\cal {R}}_{0}$ be the subring of the field of rational functions in $v$, $s$ generated by $v$, $s$, $v^{-1}$, $s^{-1}$, $\dfrac{1}{s^{2n}-1}$ for $n>0$, $\dfrac{1}{s^{2m}-v^2}$ for all $m$.
\begin{Prop}
Let $k={\cal {R}}_{0}[x,x^{-1}]$, then 
$S_{0}(S^{1} \times S^2)$ is the free $k$-module generated by $\phi$.
\end{Prop}

{\em Remark}. Proposition 1.4 allows us to define a ``Homflypt rational function'' in ${\cal {R}}_{0}[x,x^{-1}]$ for framed oriented null homologous links in $S^1\times S^2$. If $L$ is such a link, one defines $f(L)$ by $L=f(L)\phi\in S_{0}(S^{1} \times S^2)$. 
%\end{document}
%Sample file: article.tpl
%Typeset with AMSlaTex format file
%Preamble 
% Style section 

%\begin{document}
\section{A basis for $S(S^1 \times D^2, A, B)$}

\subsection{Motivation}
Suppose $M$ is an oriented $3$-manifold with boundary $\partial M$ and $\gamma$ is a nontrivial simple closed curve in $\partial M$. Let $N$ be the $3$-manifold obtained from $M$ by attaching a $2$-handle along $\gamma$. The natural inclusion $i$: ${M} \to {N}$ induces an epimorphism $i_{*}$: ${S(M)} \to {S(N)}$.  The effect of adding a $2$-handle to $M$ is to add relations to the module ${S(M)}$. Our goal is to identify {\em ker} $i_{*}$. Following G. Masbaum's work in the case of the Kauffman bracket skein module \cite{gM96}, we will use the following method to parametrize the relations coming from sliding over the $2$-handle. Pick two points $A, \ B$ on $\gamma$, which decompose $\gamma$ into two intervals $\gamma '$ and $\gamma ''$. 
%\hfil {\epsffile{gab.ai}} \hfil

Let $R$ be the submodule of $S(M)$ given by the collection $\{\Phi '( z ) - \Phi ''( z ) \mid z \in S(M, A, B)\}$. Here $z$ is any element of the relative skein module $S(M, A, B)$, and $\Phi '( z )$ and $\Phi ''( z )$ are given by capping off with $\gamma '$ and $\gamma ''$, respectively, and pushing the resulting links back into $M$.

Let $X$ be an oriented $3$-manifold. Let $F(X)$ denote the set of framed oriented links in $X$, let $I(X)$ denote the equivalence relation given by isotopy of framed oriented links in $X$; let ${\frak{B}}(X)$ denote the equivalence relation on $F(X)$ given by $L\sim bL$ if $L$ is a framed oriented link and $bL$ is obtained from $L$ by band summing $L$ with $\gamma$ at $A, B$ along some band. 

Let ${\cal{F}}(X)$ be the $k$-module generated by $F(X)$; let ${\cal{B}}(X)$ be the $k$-submodule of ${\cal{F}}(X)$ generated by elements of the form $L-bL$ where $L$ and $bL$ are as above; Let ${\cal{I}}(X)$ be the $k$-submodule of ${\cal{F}}(X)$ generated by elements of the form $L_1-L_2$, where $L_1$ and $L_2$ are isotopic. Let ${\cal{H}}(X)$ denote the submodule of ${\cal{F}}(X)$ generated by the Homflypt skein relations; let $\tilde{\cal{H}}(X)$ denote  ${\cal{H}}(X)/({\cal{I}}(X)\cap {\cal{H}}(X))$. Let $\tilde{\cal{B}}(X)$ denote the module ${\cal {B}}(X)$ modulo $({\cal{I}}(X)+{\cal{H}}(X))\cap {\cal {B}}(X)$. It is not hard to see that $\tilde{\cal{B}}(M)$ is $R$.

\begin{Thm}
$S(N)\cong  S(M)/ R$.
\end{Thm}
\begin{pf}
The embedding $i$: ${M} \to {N}$ induces two maps $F(M)\to F(N)$ and $F(M)/I(M)\twoheadrightarrow F(N)/I(N)$. The second map is surjective by a general position argument. Also by a general position argument, 
\begin{equation}
F(N)/I(N) \cong F(M)/I(M)* {\frak{B}}(M) \tag{2.1}
\end{equation}
Here $I(M)* {\frak{B}}(M)$ denotes the equivalence relation generated by $I(M)$ and ${\frak{B}}(M)$. One has that
\begin{equation}
k(F(M)/(I(M)*{\frak{B}}(M)))\cong {\cal{F}}(M)/({\cal{I}}(M)+{\cal{B}}(M))\tag{2.2}
\end{equation}
On the other hand, by definition, 
$S(N)\cong {\cal{F}}(N)/({\cal{I}}(N)+{\cal {H}}(N))$

$\cong ({\cal{F}}(N)/{\cal{I}}(N))/{\tilde{\cal {H}}}(N)$

$\cong k(F(N)/I(N))/\tilde{\cal{H}}(N)$

$\cong k(F(M)/I(M)* {\frak{B}}(M))/\tilde{\cal{H}}(N)$, this is by Equation (2.1). 

$\cong k(F(M)/I(M)* {\frak{B}}(M))/\tilde{\cal{H}}(M)$, this is by a general position argument.

$\cong ({\cal{F}}(M)/({\cal{I}}(M)+{\cal{B}}(M)))/\tilde{\cal{H}}(M)$, this is by Equation (2.2).

$\cong {\cal{F}}(M)/({\cal{I}}(M)+{\cal{H}}(M)+{\cal{B}}(M))$

$\cong ({\cal{F}}(M)/({\cal{I}}(M)+{\cal{H}}(M)))/\tilde{\cal{B}}(M)$

$\cong S(M)/\tilde{\cal{B}}(M)\cong S(M)/ R$.

Thus $S(N)\cong S(M)/R$. \end{pf}

{\em Remark}: (1) It follows that {\em ker} $i_{*}=R$. (2) Suppose $N_1$ is obtained by adding a $3$-handle to an oriented $3$-manifold $M_1$ with boundary $\partial M_1$. It is well known that the map $S(M_1) \to S(N_1)$ induced by inclusion is an isomorphism. This follows easily from the general position argument. (3) The space $S^1 \times S^2$ can be obtained from the solid torus $S^1 \times D^2$ by first attaching a $2$-handle along the meridian $\gamma$ and then attaching a $3$-handle. We have:
\begin{Cor}
$S(S^1\times S^2)\cong  S(S^1\times D^2)/ R$.
\end{Cor}
\noindent Here $R=\{\Phi '( z ) - \Phi ''( z ) \mid z \in S(S^1 \times D^2, A, B)\}$, $z$ is any element of the relative skein module $S(S^1 \times D^2, A, B)$, and $\Phi '( z )$ and $\Phi ''( z )$ are given by capping off with $\gamma '$ and $\gamma ''$, respectively, and pushing the resulting links back into $S^1 \times D^2$.

\subsection{ A basis for the relative Homflypt skein module of $S(S^1 \times D^2, A, B)$}

By our previous definition for the relative skein module, $S(S^1 \times D^2, A, B)$ is generated by isotopy classes of framed oriented links in $S^1 \times D^2$ with boundary an input point $A$ and an output point $B$. Such links consist of a collection of framed oriented closed curves and a framed oriented arc joining the two points $A$ and $B$. 

Let $i \in {\Bbb{Z}}, i >0$, and let $A_{i}'$ be the following element in $S(S^1 \times D^2, A, B)$ 
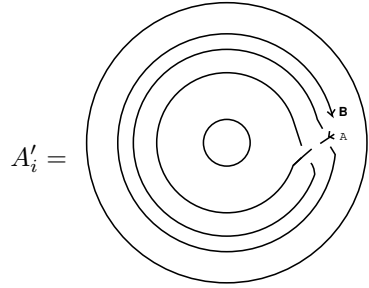
$$
A_{i}'=\raisebox{-25mm}{\ \epsfxsize1.5in\epsffile{aia.ai}}.
$$
Here if we connect $A, \ B$ by a straight line segment in $A_{i}'$, we get $A_{i}$.
Let
$$
A_{-i}'=\raisebox{-20mm}{\ \epsfxsize1.5in
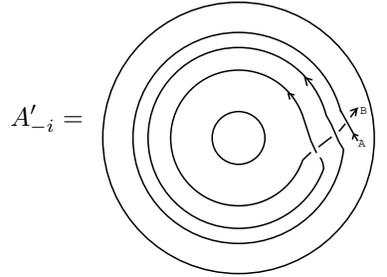
\epsffile{ai1a.ai}}.
$$
Here connecting $A, B$ by a straight line segment, we have the following:
$$
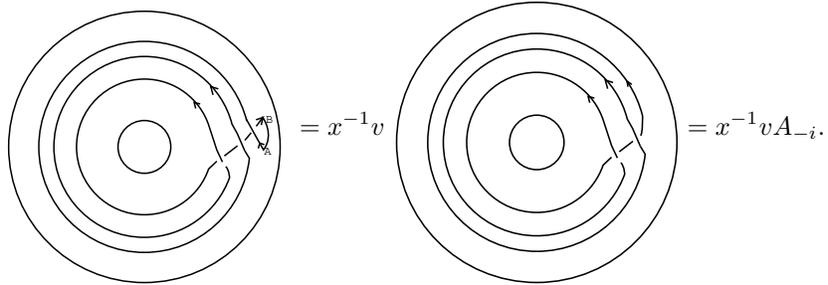
\raisebox{-20mm}{\ \epsfxsize1.5in\epsffile{ai1a1.ai}}=x^{-1}v\raisebox{-20mm}{\ \epsfxsize1.5in\epsffile{ai1a2.ai}}=x^{-1}v A_{-i}.
$$

Let $A_{0}'$ be the element given by:
$$
A_{0}'=\raisebox{-20mm}{\ \epsfxsize1.5in
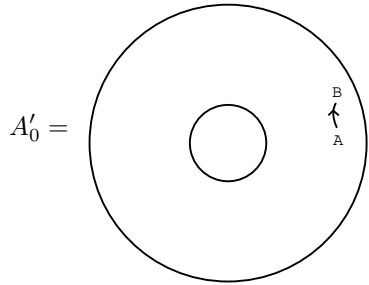
\epsffile{a0a.ai}}.
$$

\begin{Thm}
The monomials $(A_{i_1}A_{i_2}\cdots A_{i_k})A_{i}'$  for $i\in {\Bbb Z}, i_{\alpha} \in {\Bbb Z}-\{0\}, 0\leq \alpha \leq k, k \geq 0$  form a linear basis for the relative Homflypt skein module $S(S^1 \times D^2, A, B)$.
\end{Thm}
Recall $A_{i_1}A_{i_2}\cdots A_{i_k}$ is a basis element of the monomial basis of $S(S^1 \times D^2)$. We will call $(A_{i_1}A_{i_2}\cdots A_{i_k})A_{i}'$ for $i>0$ the type $1$ monomial generators, and $(A_{i_1}A_{i_2}\cdots A_{i_l})A_{-i}'$ for $i\geq 0$ the type $2$ monomial generators for $S(S^1 \times D^2, A, B)$. An example $A_{-1}A_{2}A_{2}' $ of a type 1 monomial generator is shown as:
$$
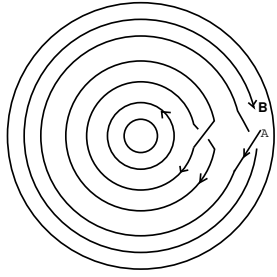
\raisebox{-25mm}{\ \epsfxsize1.5in\epsffile{ab2.ai}}
$$
 
We will use the following diagrams to illustrate the type $1$ and type $2$ monomial generators. 
$$
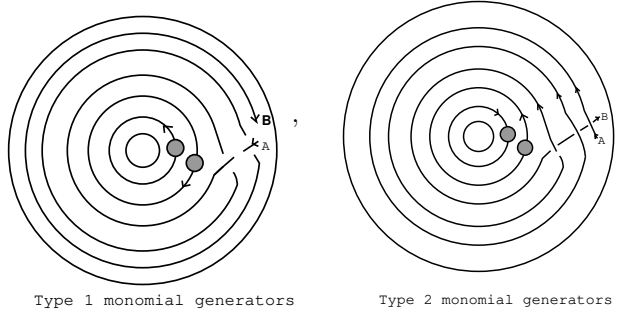
\raisebox{-34mm}{\ \epsfxsize1.5in\epsffile{ab.ai}},\quad \raisebox{-25mm}{\ \epsfxsize1.5in\epsffile{ba.ai}}
$$

Here we use a single string with a shaded circle to indicate all monomial basis elements in $S(S^1 \times D^2)$ with the given orientation.
\begin{pf}
Given a relative framed oriented link in $S(S^1 \times D^2, A, B)$, to the diagram D of this link, we can assign an ordering to the link components, where the arc is the first in the order. We also take a distinguished base point on each component, where the arc is based at $A$. The associate descending link diagram is obtained from D by changing the crossings so that when traveling around all the components in the assigned order, always beginning at the base point of each component, and each crossing is first encountered as an over-pass.

For a given relative framed oriented link, a change of  crossing can be realized by applying the Homflypt skein relations. By induction on the number of crossings, one can show that $L$ can be written as a linear combination of descending closed curves and a descending arc in $S^1 \times D^2$, where the descending closed curves lie below the arc. Any such descending link (up to change of framing) will be one of our generators. This is a sketch of a proof that $S(S^1 \times D^2, A, B)$ is generated by the type 1 and type 2 monomial generators. 

We need only to show the linear independence. Now suppose that the collection of all the monomial generators of type 1 and type 2 are linearly dependent. i.e. there is a finite linear sum
$$\sum_i \bigl(\sum_{(i_1, i_2,\dots, i_k)}a_{i_1 i_2\cdots i_k}A_{i_1}^{e_{i_1}}A_{i_2}^{e_{i_2}}\cdots A_{i_k}^{e_{i_k}}\bigr)A_{i}'=0 
$$
where all the coefficients $a_{i_1 i_2\cdots i_k}$ are nonzero.

Let $N_1$ be the largest integer appearing as a subscript to $A'$, and let  $N_2$ be the largest integer appearing as a subscript to A in the above linear combination. Let $M=\max\ \{N_2+1- N_1, 1\}$.

Now introduce a wiring of $S(S^1 \times D^2, A, B)$ into $S(S^1\times D^2)$ of the type indicated below. 
$$
\raisebox{-31mm}{\ \epsfxsize0.0in\epsffile{abtoab.ai}} 
$$
i.e. 
$$
\raisebox{-31mm}{\ \epsfxsize1.5in\epsffile{ab.ai}} \to \raisebox{-25mm}{\ \epsfxsize1.5in\epsffile{abab.ai}}
$$
$$
\raisebox{-20mm}{\ \epsfxsize1.5in\epsffile{ba.ai}} \to \raisebox{-20mm}{\ \epsfxsize1.5in\epsffile{abab1a.ai}}
$$
Such that under the wiring, $(A_{i_1}^{e_{i_1}}A_{i_2}^{e_{i_2}}\cdots A_{i_k}^{e_{i_k}})A_{i}'$  is sent to a nonzero scalar multiple of $(A_{i_1}^{e_{i_1}}A_{i_2}^{e_{i_2}}\cdots A_{i_k}^{e_{i_k}})A_{i+M}
$. In particular $A_{N_1}$ to a nonzero scalar multiple of $A_{N_1+M}$.

Now under the wiring map, the image of $f$ is
$$\sum_i \bigl(\sum_{(i_1, i_2,\dots, i_k)}a_{i_1 i_2\cdots i_k}'A_{i_1}^{e_{i_1}}A_{i_2}^{e_{i_2}}\cdots A_{i_k}^{e_{i_k}}\bigr)A_{i+M}=0 
$$
which gives a linear relation for the monomial basis elements of $S(S^1\times D^2)$, where $a_{i_1 i_2\cdots i_k}'$ is a nonzero scalar multiple of $a_{i_1 i_2\cdots i_k}$.

By our definition of $N_1, N_2$, the element $A_{N_1+M}$ will only appear in the following term.
$$
\bigl(\sum_{(N_{1_1}, N_{1_2},\dots, N_{1_k})}a_{N_{1_1} N_{1_2}\cdots N_{1_k}}'A_{N_{1_1}}^{e_{N_{1_1}}}A_{N_{1_2}}^{e_{N_{1_2}}}\cdots A_{N_{1_k}}^{e_{N_{1_k}}}\bigr)A_{N_1+M}
$$
By the linear independence of the monomial basis of $S(S^1\times D^2)$, we have 
$$a_{N_{1_1} N_{1_2}\cdots N_{1_k}}'=0,\ i.e.\ a_{N_{1_1} N_{1_2}\cdots N_{1_k}}=0$$
This contradicts the hypothesis that all the coefficients are nonzero.
\end{pf}

\subsection{Alternative type 2 monomial generators}

To get a more convenient set of relations, we find a new basis for $S(S^1 \times D^2, A, B)$.
\begin{Thm}
The following 3 types of generators form an alternative basis for $S(S^1 \times D^2, A, B)$.
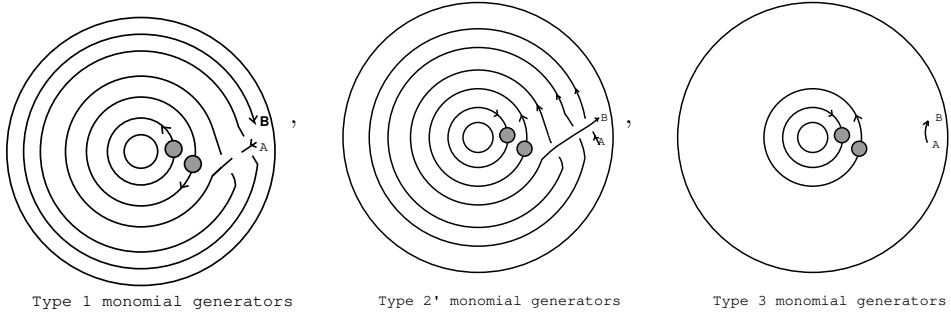
$$
\raisebox{-34mm}{\ \epsfxsize1.5in\epsffile{ab.ai}}, \quad
\raisebox{-25mm}{\ \epsfxsize1.5in\epsffile{ba1.ai}}, \quad \raisebox{-25mm}{\ \epsfxsize1.5in\epsffile{baa.ai}}
$$
\end{Thm}
Here the type 2' monomial generator has an outside arc which goes at least once around $S^1\times D^2$.
\begin{pf} We show that there is a triangular change of basis from the set of type 2 monomial generators given in Theorem 2.3 to the set of type 2' and type 3 generators listed above. 

For each type 2 generator, we can change the negative crossings involved in the outside arc into positive crossings by applying the Homflypt skein relations. Note that smoothing the first negative crossing will result in an $A_{-1}$ component and an outside arc of fewer crossings. For instance for the case $A_{-3}'$, we have
$$
A_{-3}'=\raisebox{-20mm}{\ \epsfxsize1.5in\epsffile{ai1a.ai}}
$$
$$
=x^{-2}\raisebox{-20mm}{\ \epsfxsize1.5in\epsffile{ai11.ai}}-x^{-1}(s-s^{-1})\raisebox{-20mm}{\ \epsfxsize1.5in\epsffile{ai12.ai}}
$$
\noindent where the last term is $A_{-1}A_{-2}'$. Repeating the same process to the intermediate diagrams, we can change all negative crossings and write $A_{-3}'$ as a linear combination of elements of the forms given in the theorem.

In general, by the same process, we can change all negative crossings in $A_{-i}'$ for $i>0$ and write it as a linear combination of type 2' and type 3 generators.
One can see that this is a triangular change of basis, the result follows.
\end{pf}
By Corollary 2.2, $S(S^1\times S^2)\cong  S(S^1\times D^2)/ R$, where $R=\{\Phi '( z ) - \Phi ''( z ) \mid z \in S(S^1 \times D^2, A, B)\}$. If $\alpha, \beta \in S(S^1 \times S^2)$ and $\alpha-\beta \in R$, we will say $\alpha \equiv \beta$ is a relation. We will say a set of relations $\{\alpha_i \equiv \beta_i\}$ is complete if the set $\{\alpha_i-\beta_i\}$ generates $R$. Therefore we have the following theorem.

\begin{Thm}
The following is a complete set of relations:
\begin{equation}
\raisebox{-34mm}{\ \epsfxsize1.5in\epsffile{ab1.ai}} \equiv \raisebox{-30mm}{\ \epsfxsize1.5in\epsffile{ab1a.ai}} \tag{2.3}
\end{equation}
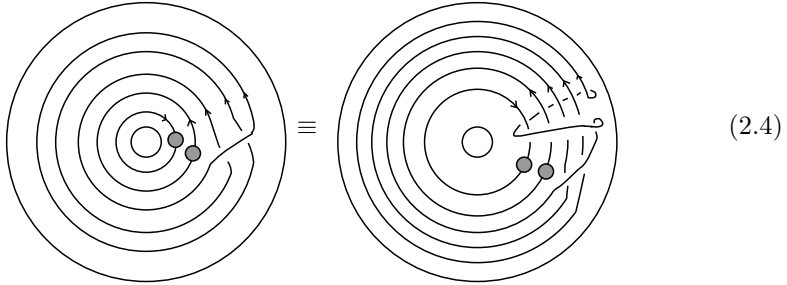
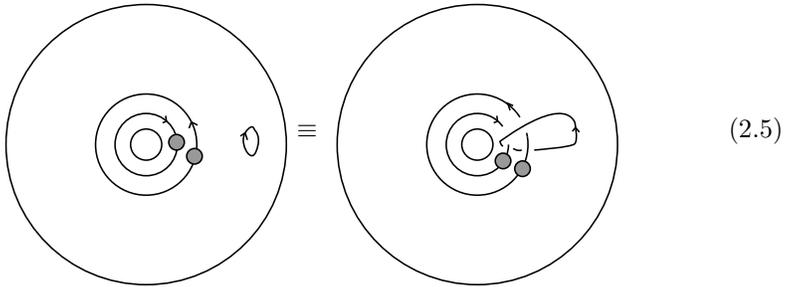
\begin{equation}
\raisebox{-20mm}{\ \epsfxsize1.5in\epsffile{ba11.ai}} \equiv \raisebox{-20mm}{\ \epsfxsize1.5in\epsffile{ba15.ai}} \tag{2.4}
\end{equation}
\begin{equation}
\raisebox{-20mm}{\ \epsfxsize1.5in\epsffile{baa1.ai}} \equiv \raisebox{-20mm}{\ \epsfxsize1.5in\epsffile{baa2.ai}} \tag{2.5}
\end{equation}
\end{Thm}
\begin{pf}
Taking $z$ to be the type 1 and type 3 monomial generators of $S(S^1 \times D^2, A, B)$, Equations (2.3) and (2.5) follow directly from the Corollary. Taking $z$ to be the type 2' monomial generators of $S(S^1 \times D^2, A, B)$, we have:
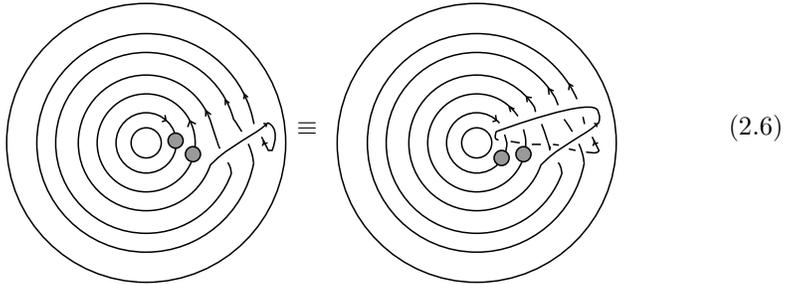
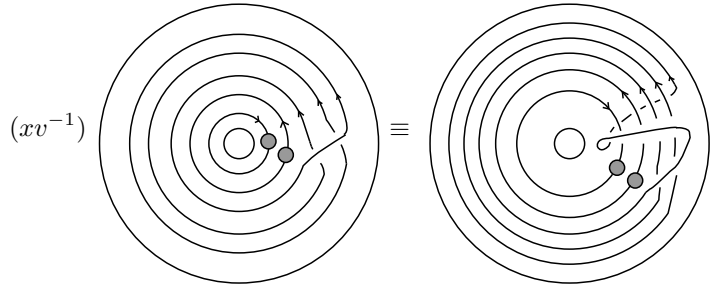
\begin{equation}
\raisebox{-20mm}{\ \epsfxsize1.5in\epsffile{ba1a.ai}} \equiv \raisebox{-20mm}{\ \epsfxsize1.5in\epsffile{ba1aa.ai}} \tag{2.6}
\end{equation}
i.e.
$$
(xv^{-1})\raisebox{-20mm}{\ \epsfxsize1.5in\epsffile{ba11.ai}} \equiv \raisebox{-20mm}{\ \epsfxsize1.5in\epsffile{ba12.ai}}
$$
Therefore
\begin{equation}
\raisebox{-20mm}{\ \epsfxsize1.5in\epsffile{ba11.ai}}\equiv (x^{-1}v)^2\raisebox{-20mm}{\ \epsfxsize1.5in\epsffile{ba13.ai}} \tag{2.7}
\end{equation}
This is equivalent to (2.4).
\end{pf}
{\em Remark:}
$$
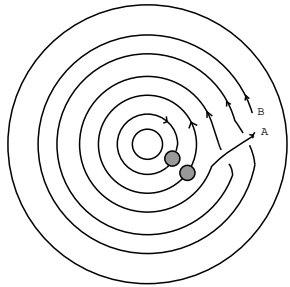
\raisebox{-20mm}{\ \epsfxsize1.5in\epsffile{ba14.ai}}
$$

Equation (2.4) identifies the wiring images of elements of the above type in $S(S^1 \times D^2, B, A)$ under two different wirings.

The relations in the above presentation are difficult to compute and analyse explicitly. For this reason, we need a new basis for $S(S^1 \times D^2, A,B)$ related to the Young idempotent basis of $S(S^1 \times D^2)$. We will use it to partially compute the relations. We remark that one can easily describe presentations for the Homflypt skein modules of lens spaces similar to that given in Theorem 7.

%\end{document}
%\begin{document}
\section{preliminaries for Young idempotents}

This section (except for Corollary 3) is a summary for the related work by C. Blanchet \cite{cB98}, A. Aiston and H. Morton \cite{A96}, \cite{AM98}. See their papers for further references to the origin of some of these ideas and results in the work of others. {\em From now on, we will assume that $s^{2n}-1$ is invertible for integers $n>0$}. It follows that the quantum integers $[n]=\dfrac{s^n-s^{-n}}{s-s^{-1}}$ for $n> 0$ are invertible in $k$. Let $[n]!=\prod_{j=1}^{n}[j]$, so $[n]!$ is invertible for $n>0$.

\subsection{Idempotents in the Hecke Algebra} 

Recall that, to a
partition of $n$, $\lambda=(\lambda_1\geq \dots \lambda_p \geq
1)$, $\lambda_{1}+\dots +\lambda_{p }=n$, is associated a Young
diagram of size $|\lambda|=n$, which we denote also by $\lambda$.
This diagram has $n$ cells indexed by \{$(i, j), \ 1\leq i\leq p$,
$1\leq j\leq \lambda_i$\}. If $c$ is the cell of index $(i, j)$ in
a Young diagram $\lambda$, its hook-length $hl(c)$ and its content
$cn(c)$ are defined by
$$hl(c)=\lambda_i+\lambda_j^{\vee}-i-j+1,\quad cn(c)=j-i.$$ Here
$\lambda^{\vee}$ is the transposed Young diagram of $\lambda$, and
$\lambda_{j}^{\vee}$ is the length of the $j$-th column of
$\lambda$ (the $j$-th row of $\lambda^{\vee}$).

For a Young diagram $\lambda$, we will use the notation $[hl(\lambda)]$, for the product over all cells of the quantum hook-lengths.
$$[hl(\lambda)]=\prod_{cells}[hl(c)]$$

\begin{Def}
{\bf The Hecke category.} The $k$-linear Hecke category $H$ is defined as follows. An object in this category is a disc $D^2$ equipped with a set of framed points. If $\alpha=(D^2, l)$ and $\beta=(D^2, l')$ are two objects, the module $Hom_H(\alpha, \beta)$ is ${\cal{S}}(D^2 \times [0,1], l\times 1, l'\times 0)$. The notation $H(\alpha, \beta)$ and $H_{\alpha}$ will be used for $Hom_H(\alpha, \beta)$  and $H(\alpha, \alpha)$ respectively. The composition of morphisms are by stacking the first one on the top of the second one.
\end{Def}
C. Blanchet stacks the second one on the top of the first. His strands are oriented upward, whereas ours are oriented downward. We arrange the rows of a Young diagram with shorter rows beneath the longer rows as is usual. This also differs from C. Blanchet's convention. Our conventions agree with those of H. Morton and A. Aiston in these regards.

Let $\otimes$ denote the monoid structure on $H$ given by embedding two disks $D^2$ side by side into one disk. For a Young diagram $\lambda$, let $\square_{\lambda}$ denote the object of the category $H$ obtained by assigning each cell of ${\lambda}$ a point equipped with the horizontal (to the left) framing. When $\lambda$ is the Young diagram with a single row of $n$ cells, $H_{\square_{\lambda}}$ will be denoted by $H_n$. $H_n$ is the $n$th Hecke algebra of type A. \cite{MT90}, \cite{vT90}.

\begin{Def}
\cite{MT90}
{\bf A positive permutation braid} is defined for each permutation
$\pi\in S_n$. It is the n-string braid, $w_{\pi}$, uniquely
determined by the properties

(1) all strings are oriented from top to bottom,

(2) for $i=1, \cdots, n$, the $i$th string joins the point
numbered $i$ at the top of the braid to the point $\pi(i)$ at the
bottom of the braid,

(3) all the crossings occur as positive crossings and each pair
crosses at most once.
\end{Def}

{\em Symmetrizers.} Let $\sigma_i \in H_n,\ i=1,\dotsc,n-1$, be the positive permutation corresponding to the transposition $(i\ i+1)$. The following theorem is shown in \cite{cB98} and \cite{AM98}.
\begin{Thm}
For $n>0$, there exist unique idempotents $f_n, \ g_n\in H_n$ such that  $\sigma_{i}f_n=xsf_n=f_n\sigma_i$ and  $\sigma_{i}g_n=-xs^{-1}g_n=g_n\sigma_i$ for all $1\leq i\leq n-1$. Moreover we have
$$f_n=\dfrac{1}{[n]!}s^{-\frac{n(n-1)}{2}}\sum_{\pi \in S_n}(xs^{-1})^{-l(\pi)}\omega_{\pi}$$
and
$$g_n=\dfrac{1}{[n]!}s^{\frac{n(n-1)}{2}}\sum_{\pi \in S_n}(-xs)^{-l(\pi)}\omega_{\pi}$$
Here $\omega_\pi$ is the positive permutation braid associated with the permutation $\pi$, and $l(\pi)$ is the length of $\pi$. \cite{M}
\end{Thm}

We summarize the Aiston-Morton description of Young symmetrizers. For a Young diagram $\lambda$ of size $n$, let $F_{\lambda}$ be the element in $H_{\square_{\lambda}}$ formed with one copy of $[\lambda_i]!f_{\lambda_i}$ along the row $i$, for $i=1,\dotsc,p$, and let $G_{\lambda}$ be the element in $H_{\square_{\lambda}}$ formed with one copy of $[\lambda_j^{\vee}]!g_{\lambda_j^{\vee}}$ along the column $j$, for $j=1,\dotsc,q$. The following four propositions are in \cite{cB98}.
\begin{Prop}
Let $\Tilde{y}_{\lambda}=F_{\lambda}G_{\lambda}$, then $\Tilde{y}_{\lambda}$ is a quasi-idempotent, and ${\Tilde{y}_{\lambda}}^2=[hl(\lambda)]\Tilde{y}_{\lambda}$. i.e. $y_{\lambda}=[hl(\lambda)]^{-1}\Tilde{y}_{\lambda}$ is an idempotent.
\end{Prop}

\begin{Prop}
Let $\lambda, \ \mu$ be two Young diagrams with $|\lambda |\ =\ |\mu|$,

(1) (Orthogonality). If $\lambda \ \ne \ \mu$, then $y_{\lambda}H({\square_{\lambda}}, {\square_{\mu}})y_{\mu}\ =\ 0$.

(2) $y_{\lambda}H_{\square_{\lambda}}y_{\lambda}=cy_{\lambda}$, where $c$ is a scalar.
\end{Prop}

\begin{Prop}
(Absorbing property). Let $\lambda \subset \mu$ be two Young diagrams, the complement of $\lambda$ in $\mu$ is called a skew Young diagram and is denoted by $\mu / \lambda$. One has
$ y_{\mu}\rho^{-1}( y_{\lambda} \otimes 1_{\square_{\mu / \lambda}})\rho y_{\mu}\ =\ y_{\mu}$.
\end{Prop}
Here $\square_{\mu / \lambda}$ is an object in $H$ with a point assigned in each cell of $\mu / \lambda$, $\rho$ is an isomorphism given by moving each point in $\mu / \lambda$ to its position in $\mu$. We will only apply this in the case $|\mu|= |\lambda | + 1$.

\begin{Prop}
(a) Let $\mu'$ be obtained by deleting an extreme cell $c$ from $\mu$. Note  $|\mu|=|\mu'|+1$. Then 
$$
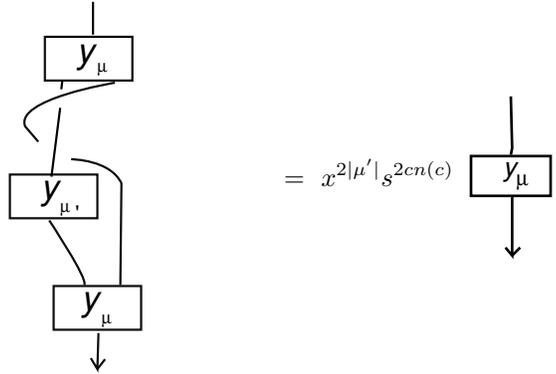
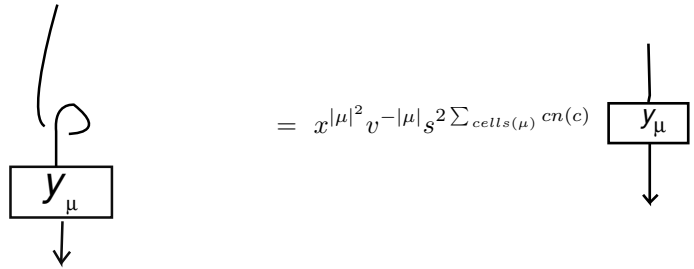
\raisebox{-25mm}{\ \epsfxsize1.5in\epsffile{fig1.ai}} =\ x^{2|\mu '|}s^{2cn(c)}\raisebox{-10mm}{\ \epsfxsize0in\epsffile{fig2.ai}}.
$$
(b) (Framing coefficient)
$$\raisebox{-18mm}{\ \epsfxsize1.5in\epsffile{fram.ai}} =\ x^{|\mu|^{2}}v^{-|\mu|}s^{2\sum_{cells (\mu)} cn(c)}\raisebox{-10mm}{\ \epsfxsize0in\epsffile{fig2.ai}}.$$

\end{Prop}

\begin{Cor}
If $\mu'$ and $\mu$ are as above, then 
$$\raisebox{-25mm}{\ \epsfxsize1.5in
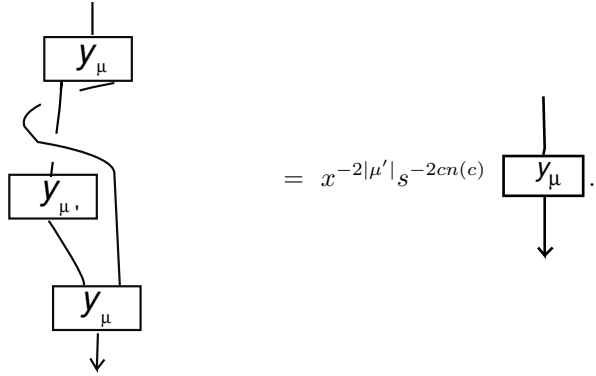
\epsffile{fig11.ai}} =\ x^{-2|\mu '|}s^{-2cn(c)}\raisebox{-10mm}{\ \epsfxsize0in\epsffile{fig2.ai}}.$$
\end{Cor}

\subsection{A Basis for the $\bf{n}$th Hecke Algebra $\bf{H_n}$}

A standard tableau $t$ with shape a Young diagram $\lambda=\lambda(t)$ is a labeling of the cells, with the integers $1$ to $n$, which is increasing along the rows and the columns. We denote by $t'$ the tableau obtained by removing the cell numbered by $n$. Note the cell numbered by $n$ in a standard tableau is always an {\em extreme cell}. C. Blanchet defines $\alpha_t \in H(n, \square_{\lambda})$ and $\beta_t \in H(\square_{\lambda}, n)$ inductively by
$$
\begin{matrix}
\alpha_1=\beta_1=1_1,\\
\alpha_t=(\alpha_{t'}\otimes 1_1)\rho_{t}y_{\lambda},\\
\beta_t=y_{\lambda}\rho_{t}^{-1}(\beta_{t'}\otimes 1_1).
\end{matrix}
$$
Here $\rho_{t} \in H(\square_{\lambda(t')}\otimes 1,\square_{\lambda})$ is the isomorphism given by an arc joining the added point to its place in $\lambda$ in the standard way. See example below.

Note that $\beta_{\tau}\alpha_t=0$ if $\tau \ne t$, and $\beta_t \alpha_t=y_{\lambda(t)}$.

\begin{Thm}
(Blanchet) The family $\alpha_t\beta_{\tau}$ for all standard tableaux $t, \tau$ such that $\lambda(t)=\lambda(\tau)$ for all Young diagrams $\lambda$ with $|\lambda|=n$ forms a basis for $H_n$.
\end{Thm}

An example of the case of $\  \lambda \ =\ \raisebox{-6mm}{\ \epsfxsize0in\epsffile{l.ai}\ }$ is as follows.
$$\hbox{Let}\  t\ =\ \raisebox{-6mm}{\ \epsfxsize0in\epsffile{zt.ai}\ }, \hbox{then}\quad   t'\ = \ \raisebox{-2mm}{\ \epsfxsize0in\epsffile{t1.ai}\ }; $$
$$\hbox{Let}\   \tau\ =\ \raisebox{-6mm}{\ \epsfxsize0in\epsffile{tau.ai}\ }, \hbox{then}\quad   \tau '\ = \ \raisebox{-6mm}{\ \epsfxsize0in\epsffile{tau1.ai}\ }; $$
The diagram description of $\alpha_{t}, \beta_{\tau}, \alpha_{t}\beta_{\tau}, \beta_{\tau}\alpha_{t}$ is the following:
$$ \alpha_{t}\ =\ \raisebox{-30mm}{\ \epsfxsize1.5in
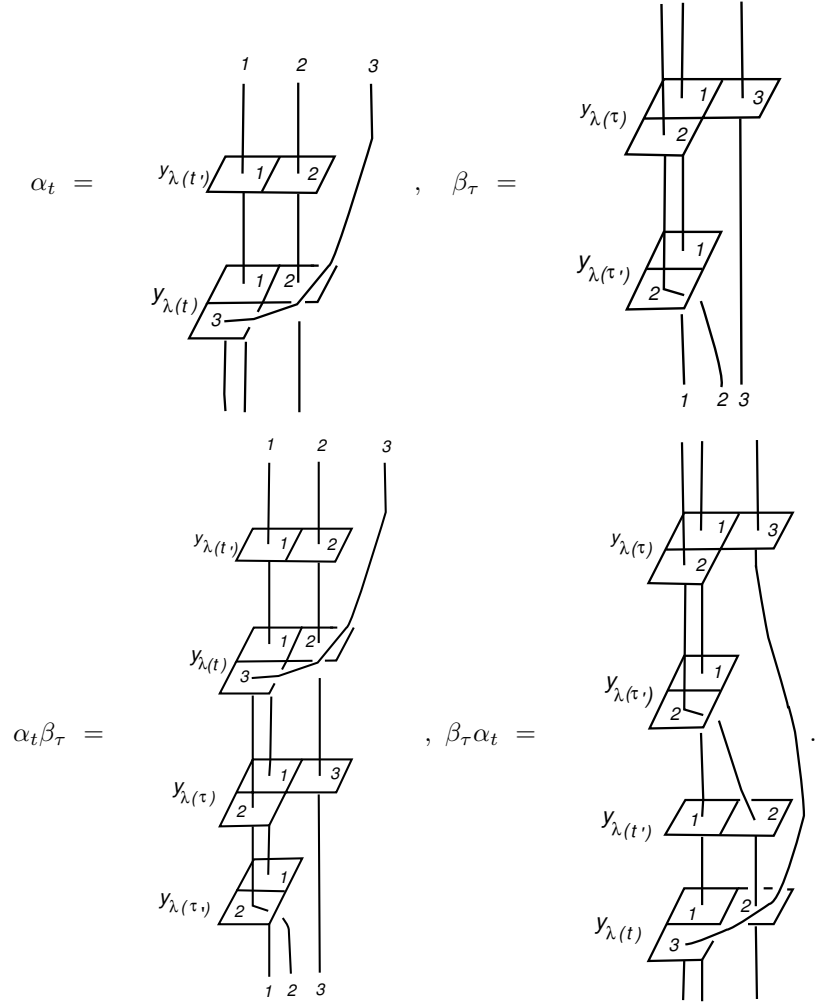
\epsffile{alphat.ai}\ }, \quad   \beta_{\tau} \ = \ \raisebox{-30mm}{\ \epsfxsize1.3in\epsffile{betau.ai}\ } $$
$$ \alpha_{t}\beta_{\tau}\ =\ \raisebox{-35mm}{\ \epsfxsize1.5in\epsffile{albet.ai}\ },\ \beta_{\tau}\alpha_{t}\ = \ \raisebox{-35mm}{\ \epsfxsize1.3in\epsffile{bealt.ai}}.$$

%Note that we obtain $\raisebox{-10mm}{\ \epsfxsize0in\epsffile{fig2a.ai}}$ %from $\raisebox{-10mm}{\ \epsfxsize0in\epsffile{fig2.ai}}$ by literally %reversing the orientation.

If $\raisebox{-10mm}{\ \epsfxsize0in\epsffile{figa.ai}}$ denotes a certain linear combination of braids, we let $\raisebox{-10mm}{\ \epsfxsize0in\epsffile{figaa.ai}}$ denote the same linear combination of braids but with the string orientation reversed on all braids which appear in the linear combination. As is usual in skein theory, a schematic diagram including several boxes which represent linear combinations of braids represents the linear combination obtained by expanding multilinearly. Let $y_{\lambda}^{*}$ denote the flattened version of $y_{\lambda}\in H_{|\lambda|}$ introduced in section 4 of \cite{AM98}. We will use $\raisebox{-10mm}{\ \epsfxsize0in\epsffile{fig22.ai}}$ to denote the flattened version of $\raisebox{-10mm}{\ \epsfxsize0in\epsffile{fig2.ai}}$. 
%\end{document}
\section{Proof of Theorem $1$ }
%\begin{document}
We introduce two types of wirings of $H_n$ into $S(S^1\times D^2, A, B)$ given by:
$$
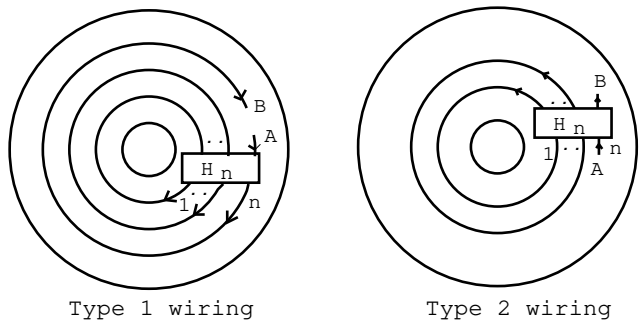
\raisebox{-25mm}{\ \epsfxsize.0in\epsffile{hnab1.ai}}\quad \raisebox{-25mm}{\ \epsfxsize.0in\epsffile{hnab2.ai}}
$$

\noindent In each case, we close the first $n-1$ strings and partially close the last string.

Let $\lambda$ be a Young diagram with $n$ cells, and let $t$ and $\tau$ be two standard tableaux of $\lambda$.  Let $\widetilde{\alpha_t\beta_{\tau}}$ be the image of $\alpha_t\beta_{\tau}$ under type 1 
wiring. Let $\overline{\alpha_t\beta_{\tau}}$ be the image of $\alpha_t\beta_{\tau}$ under type 2 
wiring.

We will try to simplify the two types of monomial generators given in chapter $2$. Note that each  type $1$ monomial generator can be written as the product of an element in $C^{-}$ and the remaining part with all its strings going in the clockwise direction. We will deal with the remaining part first, a diagram of which is shown below:
$$
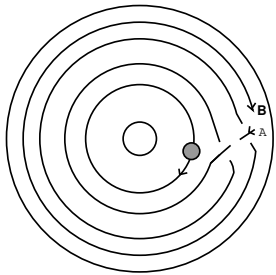
\raisebox{-25mm}{\ \epsfxsize1.5in\epsffile{abo.ai}}
$$
If we wire each element into $S(S^1 \times D^2)$ by connecting the two points $A, B$ by a straight line segment, we obtain an element in $C^{+}$.

A monomial description of these elements are $(A_{i_1}A_{i_2}\cdots A_{i_k})A_{i}'$, where the $i_1$, $i_2$, $\cdots$, $i_k$ are positive integers, $k\geq 0$ and $i$ is a positive integer. Let $C_{n}'$ denote the subspace of $S(S^1 \times D^2, A, B,)$ generated by $\{(A_{i_1}A_{i_2}\cdots A_{i_k})A_{i}'\mid i_1+ i_2+ \cdots+ i_k+i=n,\ i_1,\cdots, i_k > 0,\ i>0,\ k\geq 0\}$. Similarly, define $C_{-n}'$ to be the subspace of $S(S^1 \times D^2, A, B,)$ generated by $\{(A_{-i_1}A_{-i_2}\cdots A_{-i_k})A_{-i}'\mid i_1+ i_2+ \cdots+ i_k+i=n,\ i_1,\cdots, i_k>0,\ k\geq 0, i\geq 0\}$. A diagram description of a basis element of $C_{-n}'$ is given by:
$$
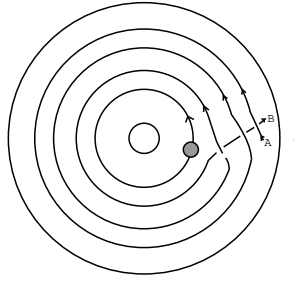
\raisebox{-18mm}{\ \epsfxsize1.5in\epsffile{bao.ai}}.
$$
\begin{Lem}
Let $n>0$, $C_{n}'$ is generated by $\widetilde{\alpha_t\beta_t}$, where $t$ takes all standard tableaux of all Young diagrams $\lambda$ with $n$ cells, and the two points $A,\ B$ are taken to be the two end points of $\widetilde{\alpha_t\beta_t}$ which have not been closed.
$$
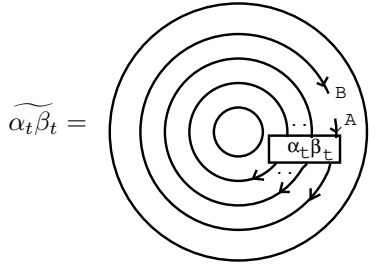
\widetilde{\alpha_t\beta_t}=\raisebox{-18mm}{\ \epsfxsize1.5in\epsffile{atbt.ai}}.
$$
\end{Lem}
\begin{pf}
First each element of $C_{n}'$ in the relative skein is isotopic to a type 1 wiring image of an $n$-string braid, which is an element in $H_n$. This element in $H_n$ can be written as a linear combination of the ${\alpha_t\beta_{\tau}}$'s, which are the generators for $H_n$. It follows from Lemma 4.3 below that ${\widetilde{\alpha_t\beta_t}}$'s are the generators for $C_{n}'$.
\end{pf}
Similarly, we have:
\begin{Lem}
If $n>0$, $C_{-n}'$ is generated by $\overline{\alpha_t\beta_t}$, where $t$ and $\lambda$ are as in the previous lemma.
\end{Lem}
$$
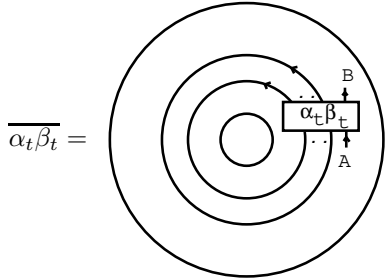
\overline{\alpha_t\beta_t}=\raisebox{-18mm}{\ \epsfxsize1.5in\epsffile{atbt1.ai}}.
$$

Since $\lambda(t)=\lambda(\tau)=\lambda$, we will write $y_{\lambda}$ for $y_{\lambda(t)}$ and $y_{\lambda(\tau)}$ in the following diagram description.

\begin{Lem}
 $\widetilde{\alpha_{t}\beta_{\tau}}=0$ if $t \ne \tau$.
\end{Lem}
\begin{pf}
 We will consider two cases:

Case $(a)$: suppose $t(n) \ne \tau(n)$, then $\lambda(t') \ne \lambda(\tau')$.
A schematic picture of $\widetilde{\alpha_{t}\beta_{\tau}}$ is as follows:
$$
{\widetilde{\alpha_t\beta_{\tau}}}\ = \ \raisebox{-20mm}{\ \epsfxsize1.5in
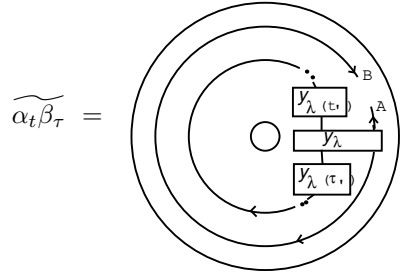
\epsffile{ylattau.ai}}
$$
The dotted part in the diagram indicates that there are other Young idempotents of Young diagrams of smaller size and arcs according to the inductive definition of $\alpha_{t}$ and $\beta_{\tau}$. We will apply the orthogonality in Proposition 3.2 to part of the picture,
$$
 \raisebox{-20mm}{\ \epsfxsize1.5in
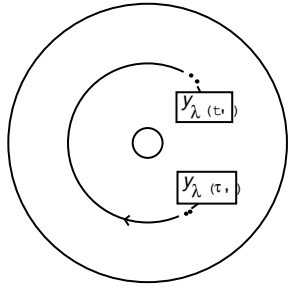
\epsffile{ylattaua.ai}},
$$
which is $0$ since $y_{\lambda(t')} H(\square _{\lambda(t')}, \square_{\lambda(\tau')}) y_{\lambda(\tau')}=0$. Therefore we have ${\widetilde{\alpha_t\beta_{\tau}}}\ = 0$.

Case $(b)$: if $t(n) = \tau(n)$, let $k$ be the integer such that $t(k+1) \ne \tau(k+1)$ and $t(k_1) = \tau(k_1)$ for $k_1 > k+1$, and let $t_k$ be the standard tableau of the Young diagram of $k$ cells obtained from $\lambda$ by removing the cells in $\lambda(t)$ labelled by $k+1, \dots, n$, we denote the obtained Young diagram by $\lambda(t_k)$. Similarly we have $\tau_k$ and $\lambda(\tau_k)$.
Since $t(k+1) \ne \tau(k+1)$, we have $\lambda(t_k) \ne \lambda(\tau_k)$.
Again we will apply the orthogonality to part of the picture of $\widetilde{\alpha_t\beta_{\tau}}$ as:
$$
 \raisebox{-20mm}{\ \epsfxsize1.6in
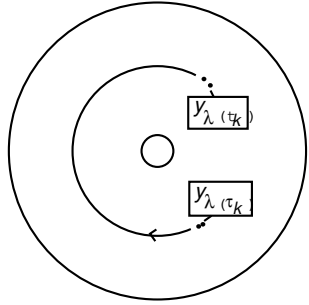
\epsffile{ylattau1.ai}}
$$
which is $0$ since $y_{\lambda(t_k)} H(\square _{\lambda(t_k)}, \square_{\lambda(\tau_k)}) y_{\lambda(\tau_k)}=0$. Therefore we have ${\widetilde{\alpha_t\beta_{\tau}}}\ = 0$.
\end{pf}
Similarly, we have:
\begin{Lem}
 $\overline{\alpha_{t}\beta_{\tau}}=0$ if $t \ne \tau$.
\end{Lem}

\begin{Lem}
 When $t = \tau$,
$$
{\widetilde{\alpha_t\beta_t}}\ = \ \raisebox{-29mm}{\ \epsfxsize1.5in\epsffile{ylatta.ai}}.
$$
\end{Lem}
\begin{pf}
Recall the inductive definition of $\alpha_t$ and $\beta_t$:
$$\alpha_t=(\alpha_{t'} \otimes 1_1)\rho_{t}y_{\lambda}$$
$$\beta_t=y_{\lambda}\rho_{t}^{-1}(\beta_{t'} \otimes 1_1)$$
Using this, a picture of the relative closure of $\alpha_t\beta_t$ is as:
$$
{\widetilde{\alpha_t\beta_t}}\ = \ \raisebox{-20mm}{\ \epsfxsize1.5in
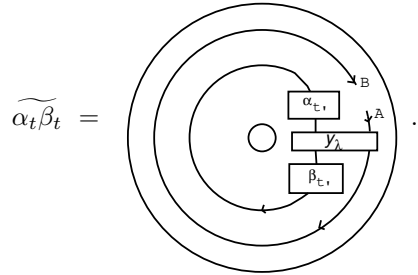
\epsffile{ylattaa.ai}}.
$$
Move $\beta_{t'}$ around to the top of $\alpha_{t'}$, applying $\beta_{t'}\alpha_{t'}=y_{\lambda(t')}$, we have
$$
{\widetilde{\alpha_t\beta_t}}\ = \ \raisebox{-20mm}{\ \epsfxsize1.5in\epsffile{ylattaaa.ai}} =\raisebox{-29mm}{\ \epsfxsize1.5in\epsffile{ylatta.ai}}.
$$
\end{pf}
A similar argument shows:
\begin{Lem}
 When $t = \tau$,
$$
{\overline{\alpha_t\beta_t}}\ = \ \raisebox{-18mm}{\ \epsfxsize1.5in\epsffile{atbt2.ai}}
$$
\end{Lem}

\begin{Prop}
Let $\widehat{\alpha_{t}\beta_{t}}$ be the closure image of $\alpha_{t}\beta_{\tau}$ of the wiring map: $H_n \to S(S^1 \times D^2)$, and $t, \ \tau$ are two standard tableaux of $\lambda$, then
$$ \widehat{\alpha_{t}\beta_{t}}=\widehat{y_{\lambda}}=\widehat{\alpha_{\tau}\beta_{\tau}}.
$$
\end{Prop}
\begin{pf}
$\widehat{\alpha_{t}\beta_{t}}$ is the closure of the string numbered by $n$ of $\widetilde{\alpha_{t}\beta_{t}}$ in $S^1 \times D^2$.  We have
$$\widehat{\alpha_{t}\beta_{t}}=\raisebox{-20mm}{\ \epsfxsize1.6in
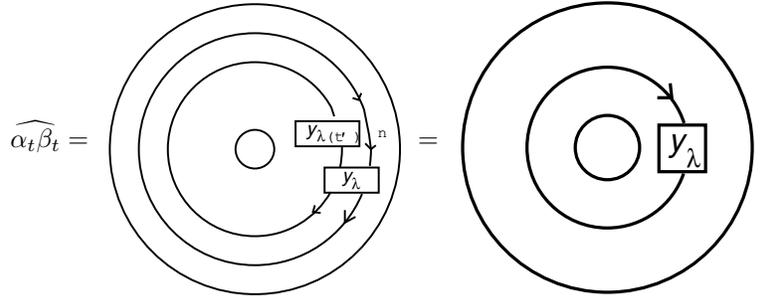
\epsffile{ylatt4.ai}}=\raisebox{-20mm}{\ \epsfxsize1.6in\epsffile{yla1.ai}}.
$$
To see this equality: (1) replace $y_{\lambda}$ by two copies of $y_{\lambda}$ using the idempotent property; (2) move the lower copy of  $y_{\lambda}$ around to the top of $y_{\lambda(t')}$; (3) use the absorbing property in Proposition 3.3. We will use this trick frequently below as simply the absorbing property. The last term is denoted by $\widehat{y_{\lambda}}$.
\end{pf}

As the closure of $H_n$ is $C_n$, one may obtain the Morton-Aiston result mentioned in the introduction as a corollary of the above proposition. Note this result is shown in \cite{A96} in a slightly different way. We used this to prove Proposition 1.1. Now put the two parts of type $1$ monomial generators together, we have

\begin{Lem}
Each type $1$ monomial generator in the relative skein module can be written as a linear sum in terms of elements of the form
$$
Q_{\lambda, \mu,c}'=\raisebox{-20mm}{\ \epsfxsize1.5in
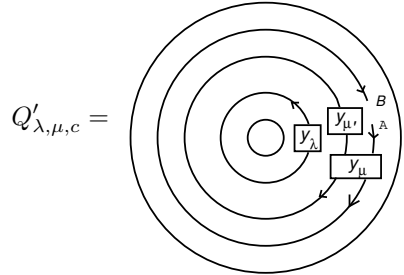
\epsffile{ylamu1.ai}}
$$
\end{Lem}

\noindent Here $y_{\mu '}$ is the Young idempotent corresponding to the Young diagram $\mu '$ obtained from the Young diagram $\mu$ by removing the extreme cell c. We will call elements of the set $\{Q_{\lambda, \mu,c}'\}$ the type 1 Young idempotent generators. In particular, if $|\lambda|=0$, we denote $Q_{\lambda, \mu,c}'$ by $Q_{0, \mu,c}'$.
\begin{pf}
Now the first part of a type $1$ monomial generator is an element in $C^{-}$, it can be written as a linear combination of generators of $C^{-}$ of the form

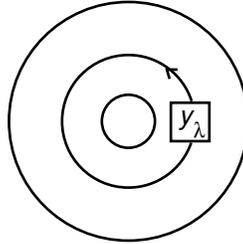
\hfil {\epsffile{yla.ai}} \hfil

It is denoted by $\widehat{y_{-\lambda}}$. When we put the two parts of type $1$ monomial generators together, we can write the type $1$ monomial generator as a linear combination of the new generators $Q_{\lambda, \mu,c}'$.
\end{pf}

Similarly, we have the following lemma.
\begin{Lem}
Each type $2$ monomial generator can be written as a linear combination of new generators of the following form.
$$
Q_{\lambda, \mu,c}''=\raisebox{-20mm}{\ \epsfxsize1.5in
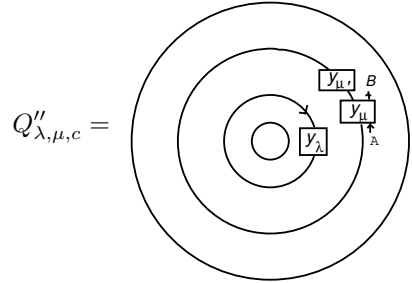
\epsffile{ylamu3.ai}}
$$
\end{Lem}
\noindent Again $y_{\mu '}$ is the Young idempotent corresponding to the Young diagram $\mu '$ obtained from the Young diagram $\mu$ by removing an extreme cell c. We will call these the type 2 Young idempotent generators for $S(S^1 \times D^2, A, B)$.

\noindent {\bf Restatement of Theorem 1}
 $S(S^1 \times D^2, A, B)$ has a basis given by \{$Q_{\lambda, \mu, c}'$, $Q_{\lambda, \mu, c}'' $\}.
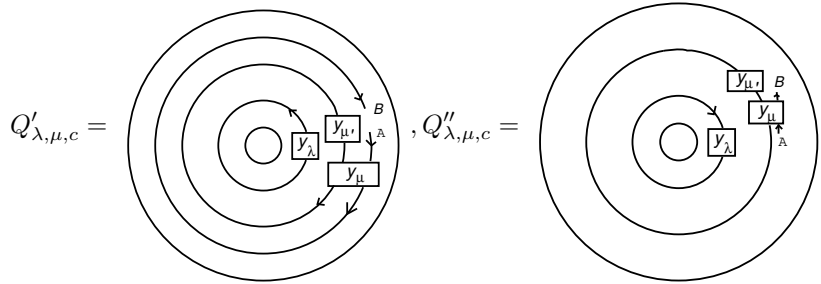
$$
Q_{\lambda, \mu, c}'=\raisebox{-20mm}{\ \epsfxsize1.5in\epsffile{ylamu1.ai}},
Q_{\lambda, \mu, c}''=\raisebox{-20mm}{\ \epsfxsize1.5in\epsffile{ylamu3.ai}}
$$

Here we change the pictures of $Q_{\lambda, \mu, c}', Q_{\lambda, \mu, c}''$ through an obvious homeomorphism of $S^1 \times D^2$.

\begin{pf}
We have shown that the type 1 and type 2 monomial generators give a basis for $S(S^1 \times D^2, A, B)$; also $S(S^1 \times D^2, A, B)$ is generated by $\{Q_{\lambda, \mu, c}'\}$ and  $\{Q_{\lambda, \mu, c}'' \}$. Moreover, each type $1$ monomial generator can be written as a linear combination of the type 1 Young idempotent generators $Q_{\lambda, \mu, c}'$ and each type $2$ monomial generator can be written as a linear combination of the type 2 Young idempotent generators $Q_{\lambda, \mu, c}''$. We will show the linear independence by comparing the dimensions.

Recall the type $1$ and type $2$ monomial generators for $S(S^1 \times D^2, A, B)$.
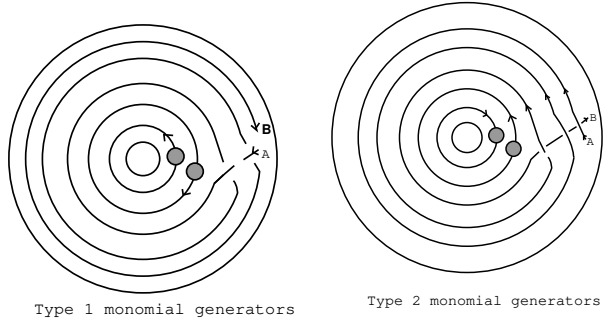
$$
\raisebox{-25mm}{\ \epsfxsize1.5in\epsffile{ab.ai}} \quad \raisebox{-15mm}{\ \epsfxsize1.5in\epsffile{ba.ai}}.
$$
Let $C_{I}$ be the subspace of $S(S^1 \times D^2, A, B)$ spanned by the type 1 monomial generators, and let $C_{II}$ be the subspace of $S(S^1 \times D^2, A, B)$ spanned by the type 2 monomial generators. We have the following:
$$C_{I}\cong \bigoplus_{m\geq 0, n>0}C_{-m}\otimes C_{n}',\quad C_{II}\cong \bigoplus_{m\geq 0, n \geq 0}C_{m}\otimes C_{-n}'.$$
Let $e(n)$ be the number of extreme cells of all the Young diagrams of size $n$, and $p(n)$ be the number of Young diagrams of size $n$, note $p(n)=\dim(C_{n})$. On one hand, we have
$$\dim(C_{n}')=\sum_{i=0}^{n-1} \dim(C_{i})=\sum_{i=0}^{n-1}p(i).$$

On the other hand, it follows from the next lemma that $$\sum_{i=0}^{n-1}p(i)=e(n).$$
Thus $e(n)=\dim(C_{n}')$.

Note $e(n)$ is also the number of diagrams of the form:
$$
Q_{0, \mu, c}'= \raisebox{-18mm}{\ \epsfxsize1.5in
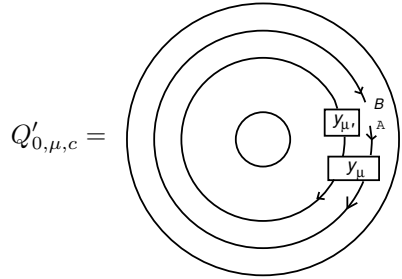
\epsffile{ylamu2.ai}}
$$
where $\mu$ is a Young diagram of size $n$ and $\mu'$ is obtained from $\mu$ by removing the extreme cell $c$.

We have shown in the proof of lemma 4.8 that $\{Q_{0, \mu, c}'\mid |\mu|=n\}$ generates $C_{n}'$. This set has cardinality $e(n)=\dim(C_{n}')$. We conclude that $C_{n}'$ has a new basis $\{Q_{0, \mu, c}'\mid |\mu|=n\}$. Recall that $C_{-m}$ has a basis given by $\{\widehat{y_{-\lambda}}\mid |\lambda|=m\}$. Hence the subspace $C_{-m} \otimes C_{n}'$ has a basis $\{Q_{\lambda, \mu, c}'\mid |\lambda|=m, |\mu|=n, {\hbox{c is an extreme cell of }} \mu\}$. Since
$$C_{I}\cong \bigoplus_{m\geq 0, n>0}C_{-m}\otimes C_{n}',$$
$C_{I}$ has a basis given by $\{Q_{\lambda, \mu, c}'\mid |\lambda|\geq 0, |\mu|>0, {\hbox{c is an extreme cell of}}\  \mu\}$.

One can study the type 2 monomial generators and the new generators $Q_{\lambda, \mu,c}''$ in a similar way. One can show that $C_{II}$ has a new basis given by $\{Q_{\lambda, \mu, c}''\mid |\lambda|\geq 0, |\mu|>0, {\hbox{c is an extreme cell of}}\  \mu\}$.

Since $S(S^1 \times D^2, A, B)\cong C_{I}\oplus C_{II} $ , the result follows.
\end{pf}

\begin{Lem}
$$e(n)=\sum_{i=0}^{n-1} p(i).$$
\end{Lem}
\begin{pf}
Consider the directed graph of the Young diagrams. The vertices are the Young diagrams, there is an edge from a Young diagram $\lambda$ to a Young diagram $\mu$ if $\mu$ can be obtained from $\lambda$ by adding a cell (which is an extreme cell for $\mu$).
$$
 \raisebox{-20mm}{\ \epsfxsize0.0in\epsffile{ygraph.ai}}
$$
For each Young diagram $\lambda$, let $i(\lambda)$ be the number of incoming edges, which is the number of ways to remove a cell from $\lambda$ to get a legitimate Young diagram. Note the removed cell is an extreme cell of $\lambda$. i.e $i(\lambda)$ is also the number of extreme cells of $\lambda$. Moreover, $i(\lambda)$ is also equal to the number of distinct lengths of rows appearing in $\lambda$. Let $o(\lambda)$ be the number of outgoing edges from $\lambda$, which is the number of ways to add a cell to $\lambda$ to get a legitimate Young diagram. There are $i(\lambda)+1$ ways to add it. So $o(\lambda)=i(\lambda)+1$. Now
$$e(n-1)=\sum_{\lambda: |\lambda|=n-1}i(\lambda)$$
And
$$e(n)=\sum_{\lambda: |\lambda|=n}i(\lambda)=\sum_{\lambda: |\lambda|=n-1}o(\lambda)=\sum_{\lambda: |\lambda|=n-1}(i(\lambda)+1)$$
Therefore we have the recursive formula for $e(n)$,

$$e(n)=e(n-1)+p(n-1).$$

Now $p(0)=1$ and $e(0)=0$, solving the recursive relation, we have
$$e(n)=\sum_{i=0}^{n-1} p(i).$$
\end{pf}

We will need the following for the next section.

If we reverse the orientations of all components of type 1 monomial generators of $S(S^1\times D^2, A, B)$ (we will refer them as the type 1 monomial generator with the reversed orientation), we get the corresponding type 1 monomial generators for $S(S^1\times D^2, B, A)$.
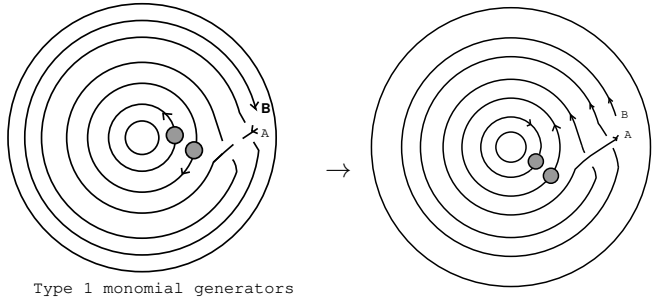
$$
\raisebox{-25mm}{\ \epsfxsize1.5in\epsffile{ab.ai}}\quad \to \raisebox{-15mm}{\ \epsfxsize1.5in\epsffile{ba14.ai}}
$$

\begin{Cor}
The subspace of $S(S^1\times D^2, B, A)$ generated by the type 1 monomial generator with the reversed orientation has an alternative basis given by

$\{\overline{Q}_{\lambda, \mu,c}'\ |\ \lambda, \mu \ {\hbox{are Young diagrams and c is an extreme cell of}}\ \mu\}$.
\end{Cor}
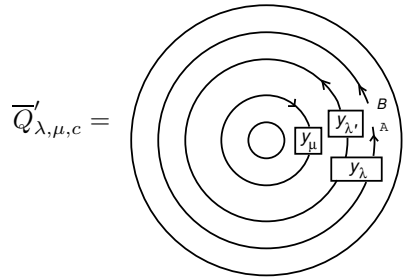
$$
\overline{Q}_{\lambda, \mu,c}'=\raisebox{-20mm}{\ \epsfxsize1.5in\epsffile{ylamu1re.ai}}
$$
Recall the relations (2.4) given in Theorem 2.7 are obtained by identifying two different wirings of generators for this subspace.
%\end{document}

%\begin{document}
\section{ Relations for the Homflypt skein module of $S(S^1 \times S^2)$ }

Recall that, $S(S^1\times S^2)\cong  S(S^1\times D^2)/ R$, where $R=\{\Phi '( z ) - \Phi ''( z ) \mid z \in S(S^1 \times D^2, A, B)\}$. It suffices to take $z$ to be the generators of $S(S^1 \times D^2, A, B)$.

We take $z$ to be the Young idempotent basis elements in $S(S^1\times D^2, A, B)$, then we have the following theorem.

\begin{Thm}
The following is a complete set of relations for $S(S^1 \times S^2)$, where $\lambda,\ \mu$ are any Young diagrams and $\mu'$ is obtained by deleting an extreme cell of $\mu$. 
\begin{equation}
\raisebox{-25mm}{\ \epsfxsize1.8in\epsffile{ylamu11.ai}}\ \equiv \ \raisebox{-25mm}{\ \epsfxsize1.8in\epsffile{ylamu12.ai}} \tag{5.1}
\end{equation}
\begin{equation}
\raisebox{-25mm}{\ \epsfxsize1.8in\epsffile{ylamu1re1a.ai}}\ \equiv \ \raisebox{-25mm}{\ \epsfxsize1.8in\epsffile{ylamu1re2a.ai}} \tag{5.3}
\end{equation}
\begin{equation}
\raisebox{-25mm}{\ \epsfxsize1.8in\epsffile{ylamua.ai}}\ \equiv \ \raisebox{-25mm}{\ \epsfxsize1.8in\epsffile{ylamuaa.ai}}\tag{5.4}
\end{equation}

\end{Thm}

{\em Remark}: Note the rest of the results of this section as well as Theorem 2 do not depend on the completeness of the above sets of relations. These results only depend on these relations being true in $S(S^1\times S^2)$, which is easily seen by sliding a strand over the $2$-handle.

Consider Equation (5.1), note that the left-hand side is equal to our generator $Q_{\lambda, \mu}$ by the absorbing property in Proposition 3.3. We will simplify the right-hand side: 
$$
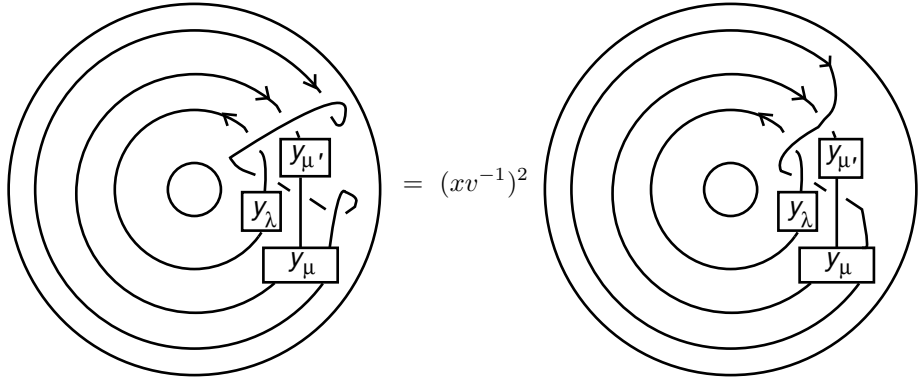
\raisebox{-25mm}{\ \epsfxsize0in\epsffile{ylamu12.ai}}\ = \ (xv^{-1})^2 \raisebox{-25mm}{\ \epsfxsize0in\epsffile{ylamu12a.ai}}
$$

\begin{Lem}
If $|\mu|\ne 0$, then $(1-x^{2|\mu|}s^{2cn(c)}v^{-2})Q_{0, \mu}\equiv 0$  for all extreme cells $c$ of $\mu$ in $S(S^1 \times S^2)$. 
\end{Lem}
\noindent Here $Q_{0, \mu}$ denotes $Q_{\lambda, \mu}$ for the case $|\lambda|=0$. Note $Q_{0, \mu}=\widehat{y_{\mu}}$.
\begin{pf}
From (5.1) with $|\lambda|=0$, we have:
$$
{\raisebox{-20mm}{\ \epsfxsize0in
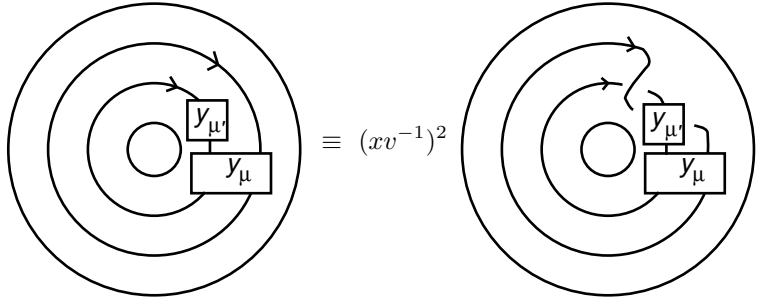
\epsffile{ymm.ai}} \ }\equiv {\ (xv^{-1})^2{\raisebox{-20mm}{\ \epsfxsize0in\epsffile{ymma.ai}}}}
$$ 
Here $c$ is the extreme cell of $\mu$ such that if we remove it, we obtain $\mu'$. The first term is equal to $Q_{0, \mu}$ by the absorbing property, the second term is equal to $(xv^{-1})^{2}x^{2|\mu '|}s^{2cn(c)} Q_{0, \mu}$ by Proposition 3.4. Therefore 
$$Q_{0, \mu}-(xv^{-1})^{2}x^{2|\mu '|}s^{2cn(c)}Q_{0, \mu} =(1-x^{2|\mu|}s^{2cn(c)}v^{-2})Q_{0, \mu}\equiv 0.$$
i.e. $Q_{0, \mu}$ is torsion.
\end{pf}

\begin{Cor}
If $\mu$ is not a rectangular Young diagram, then $Q_{0, \mu}\equiv 0$  in $S(S^1\times S^2)$.
\end{Cor}
\begin{pf}
When $\mu$ is not a rectangular Young diagram, it has at least two extreme cells $c$ and $c'$. Moreover the contents $cn(c)\ne cn(c')$. From the above lemma, 

$(1-x^{2|\mu|}s^{2cn(c)}v^{-2})Q_{0, \mu}\equiv 0$ and $(1-x^{2|\mu|}s^{2cn(c')}v^{-2})Q_{0, \mu}\equiv 0$. Therefore,

$(x^{2|\mu|}s^{2cn(c)}v^{-2}-x^{2|\mu|}s^{2cn(c')}v^{-2})Q_{0, \mu}\equiv 0$. Since $x, s$ and $v$ are invertible, we have $(s^{2cn(c)}-s^{2cn(c')})Q_{0, \mu}\equiv 0$. i.e.  $s^{2cn(c)}(1-s^{2(cn(c)-cn(c'))})Q_{0, \mu}\equiv 0$. i.e. $(1-s^{2m})Q_{0, \mu}\equiv 0$, where $m=2|cn(c)-cn(c')|$. So $Q_{0, \mu}\equiv 0$.
\end{pf}

Let $\varGamma(\lambda, \mu)=\{(\nu, \sigma) |\ |\lambda|>|\nu|,\  |\mu|>|\sigma|\  {\hbox {and}} \ |\lambda|-|\mu|=|\nu|-|\sigma|\}$.

\begin{Lem}
When $|\mu| \geq 1$, $(1-x^{2(|\mu|-|\lambda|)}v^{-2}s^{2cn(c)})Q_{\lambda, \mu}\equiv \sum_{(\nu, \sigma) \in \varGamma(\lambda, \mu)}a_{\nu\sigma}Q_{\nu , \sigma}$, where $a_{\nu\sigma} \in k$. Thus if $\mu$ is not a rectangular Young diagram, $Q_{\lambda, \mu}$ can be written as a linear combination of the $Q_{\nu, \sigma}$'s with $(\nu, \sigma) \in \varGamma(\lambda, \mu)$.
\end{Lem}
\begin{pf}
Lemma 5.1 is the case for $|\lambda|=0$. So we assume $|\lambda|\geq 1$.

By the absorbing property, the left hand side of (5.1) is $Q_{\lambda, \mu}$; the right hand side of (5.1) is $(xv^{-1})^2$ times the diagram below.
$$
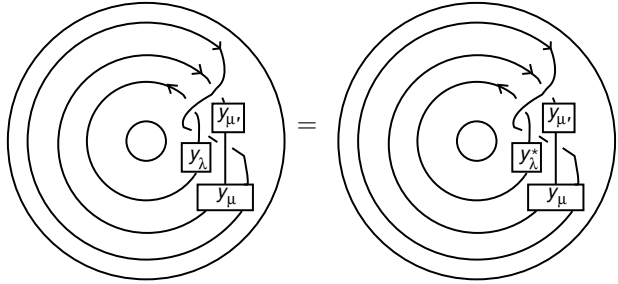
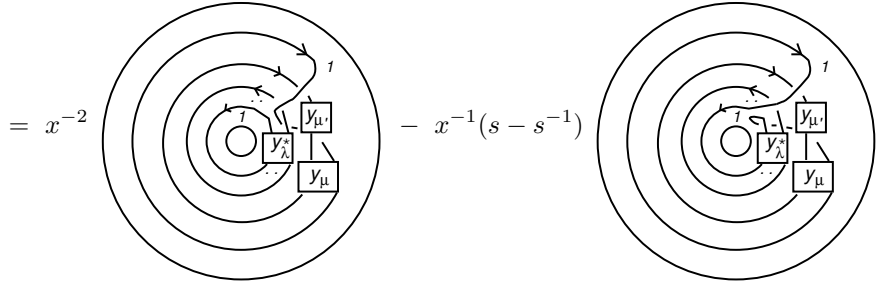
\raisebox{-20mm}{\ \epsfxsize1.5in\epsffile{ylamu12a.ai}}=\raisebox{-20mm}{\ \epsfxsize1.5in\epsffile{ylamu12aa.ai}}
$$
(As $\widehat{y_{\lambda}^{*}}=\widehat{y_{\lambda}}$),
$$=\ x^{-2}\raisebox{-20mm}{\ \epsfxsize1.5in\epsffile{y1lamua.ai}}\ - \ x^{-1}(s-s^{-1})\raisebox{-20mm}{\ \epsfxsize1.5in\epsffile{y1lamuaa.ai}}
$$
(Keep applying the skein relation to the first term on the right-hand side of the above), 
$$
= \ x^{-2}\bigl (x^{-2}\raisebox{-20mm}{\ \epsfxsize1.5in
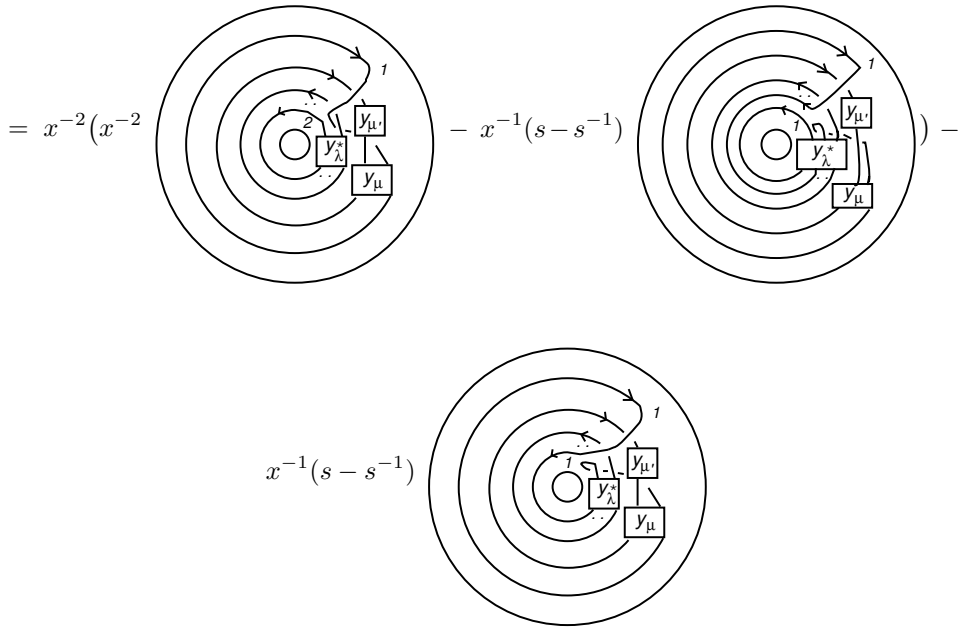
\epsffile{y2lamu.ai}}\ - \ x^{-1}(s-s^{-1})\raisebox{-20mm}{\ \epsfxsize1.5in\epsffile{y2lamua.ai}}\bigr )\ - 
$$
$$
\ x^{-1}(s-s^{-1})\raisebox{-20mm}{\ \epsfxsize1.5in\epsffile{y1lamuaa.ai}}
$$
(Repeating the above process $|\lambda|-2$ times,)
$$
= \ x^{-2|\lambda|}\raisebox{-20mm}{\ \epsfxsize1.5in\epsffile{ylamma.ai}}\ - \ \sum_{i=0}^{|\lambda|-1}x^{-2i-1}(s-s^{-1})\raisebox{-20mm}{\ \epsfxsize1.5in\epsffile{yilamu.ai}}\  
$$
(By Proposition 3.4,)
$$=(x^{-2|\lambda|}x^{2|\mu|-2}s^{2cn(c)})\raisebox{-20mm}{\ \epsfxsize1.5in\epsffile{ylamu11.ai}}\ - \ \sum_{i=0}^{|\lambda|-1}x^{-2i-1}(s-s^{-1})\raisebox{-20mm}{\ \epsfxsize1.5in\epsffile{yilamu.ai}}\  
$$
Thus
$$(1-x^{2(|\mu|-|\lambda|)}v^{-2}s^{2cn(c)})\raisebox{-20mm}{\ \epsfxsize1.5in\epsffile{ylamu11.ai}}\equiv  \sum_{i=0}^{|\lambda|-1}x^{-2i+1}v^{-2}(s-s^{-1})\raisebox{-20mm}{\ \epsfxsize1.5in\epsffile{yilamu.ai}}  
$$
Lemma 5.3 below shows that elements of the form 
$$
\raisebox{-20mm}{\ \epsfxsize1.5in\epsffile{yilamu.ai}}\ =  \raisebox{-20mm}{\ \epsfxsize1.5in\epsffile{yilamua.ai}}\ = \ x^{-1}v\raisebox{-20mm}{\ \epsfxsize1.5in\epsffile{yilamuaa.ai}}\
$$
can be written as linear combination of $Q_{\nu, \sigma}$ with $(\nu, \sigma) \in \varGamma(\lambda, \mu)$. Then the result follows.
\end{pf}

{\em Remark:} If $|\lambda|=|\mu|$, then $(1-v^{-2}s^{2cn(c)})Q_{\lambda, \mu}=\sum_{(\nu, \sigma) \in \varGamma(\lambda, \mu)}a_{\nu\sigma}Q_{\nu, \sigma}$.

\begin{Lem}
The elements of the form  
$$
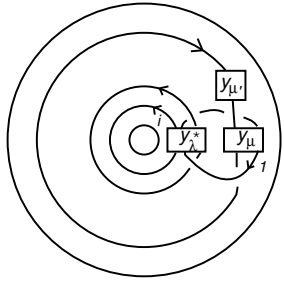
\raisebox{-20mm}{\ \epsfxsize1.5in\epsffile{yilamuaa.ai}}\  
$$
can be written as a linear combination of the $Q_{\nu, \sigma}$'s with $(\nu, \sigma) \in \varGamma(\lambda, \mu)$.
\end{Lem}
\noindent Here in the diagram $i$ is an integer such that $0\leq i \leq |\lambda|-1$.
\begin{pf}
Let $m=|\lambda|$ and $n=|\mu|$. Fix $i$, we will label two points $C, \ D$ on the diagram and label the diagram by $\lambda\mu_{CD}$. These two points separate the diagram into two parts $\lambda_{CD}$ and $\mu_{CD}$. 
$$
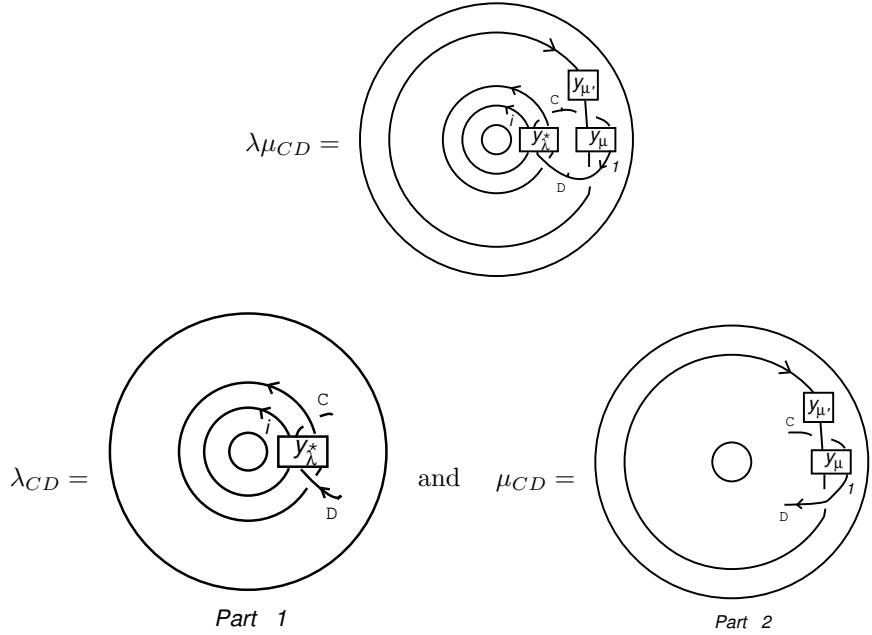
\lambda\mu_{CD}=\raisebox{-17mm}{\ \epsfxsize1.5in\epsffile{yilamucd.ai}}
$$  
$$
\lambda_{CD}=\raisebox{-20mm}{\ \epsfxsize1.5in\epsffile{yilamucd1.ai}}\quad  \hbox{and } \quad \mu_{CD}=\raisebox{-20mm}{\ \epsfxsize1.5in\epsffile{yilamucd2.ai}}
$$ 
We will consider $S^1 \times D^2$ as $X_1 \cup X_2$, where each $X_i$ is a solid torus with $C, D$ on the boundary. Note that $\lambda_{CD}\in S(X_1, D, C)$ and $\mu_{CD}\in S(X_2, C, D)$.

We will show that when we connect $\lambda_{CD}$ and $\mu_{CD}$ through $C$ and $D$, after simplification, we can rewrite $\lambda\mu_{CD}$ in terms of $Q_{\nu, \sigma}$ with $(\nu, \sigma) \in \varGamma(\lambda, \mu)$.

In particular, $\lambda_{CD}\in C_{-(m-1)}'(X_1)$, $\mu_{CD} \in C_{n-1}'(X_2)$, where $C_{-(m-1)}'(X_1)$ and $C_{n-1}'(X_2)$ are the images of $C_{-(m-1)}'(S^1\times D^2, B, A)$ and $C_{n-1}'(S^1\times D^2, A, B)$ under maps induced by obvious diffeomorphisms from $(S^1\times D^2, B, A)$ to $(X_1, D, C)$ and from $(S^1\times D^2, A, B)$ to $(X_2, C, D)$. We can use our previous study of the relative skein module with $2$ points in the boundary. It follows that $\lambda_{CD}$ can be written as a linear combination of monomial generators of $C_{-(m-1)}'(X_1)$, where the generators are given the counterclockwise orientation. The diagram descriptions of such generators are:
$$
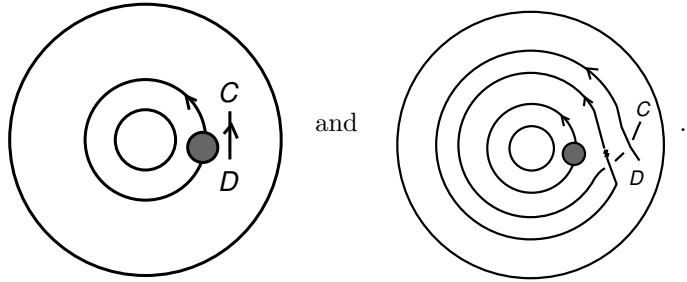
\raisebox{-20mm}{\ \epsfxsize1.5in\epsffile{type11.ai}} \quad \hbox{and} \quad \raisebox{-20mm}{\ \epsfxsize1.5in\epsffile{type31.ai}}.
$$
Similarly, $\mu_{CD}$ can be written as a linear combination of monomial generators of $C_{n-1}'(X_2)$, where the generators are given the clockwise orientation. The diagram descriptions of such generators are:
$$
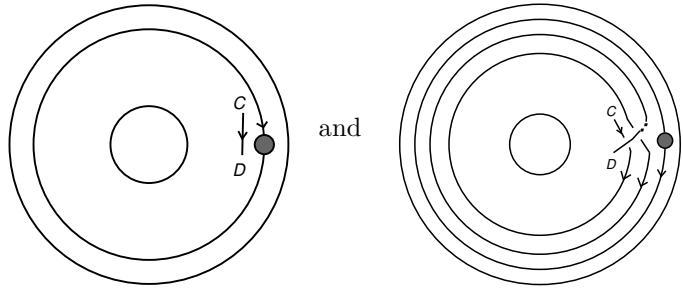
\raisebox{-20mm}{\ \epsfxsize1.5in\epsffile{type32.ai}} \quad \hbox{and} \quad \raisebox{-20mm}{\ \epsfxsize1.5in\epsffile{type42.ai}}
$$
Here again single strings with shaded circles are the monomial basis elements in $S(S^1 \times D^2)$ with the given orientation.

Therefore by the connection through $C, \ D$, $\lambda\mu_{CD}$ can be written as a linear combination of elements of the following four forms:
$$
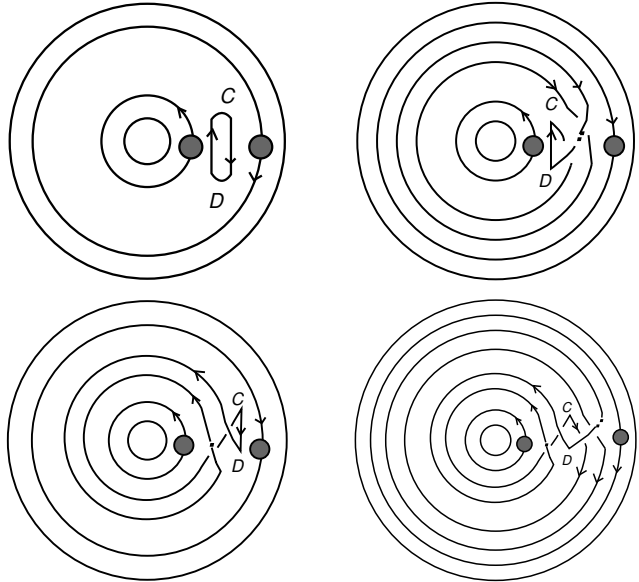
\raisebox{-20mm}{\ \epsfxsize1.5in\epsffile{type1.ai}} \quad  \quad \raisebox{-20mm}{\ \epsfxsize1.5in\epsffile{type2.ai}}
$$
$$
\raisebox{-20mm}{\ \epsfxsize1.5in\epsffile{type3.ai}} \quad  \quad \raisebox{-20mm}{\ \epsfxsize1.5in\epsffile{type4.ai}}
$$
By applying the Homflypt skein relations to the above four types of elements, the first type contains a trivial component which will contribute a scalar; and the second and the third types each contains a curl which will contribute a scalar also. Each element of the first three types is in the subspace $C_{-(m-1)} \times C_{n-1}$. For the last type, the link component containing $C, D$ is isotopic to one of the $A_i$'s, each element of the fourth type is an element in the subspaces $C_{-m_{1}} \times C_{n_1}$, where $|m_{1}|<m$ and $n_1 < n$.

Therefore each of the four types of elements can be written as a linear combination  of elements of $C_{-m_{1}} \times C_{n_1}$, where $|m_{1}|<m$ and $n_1 < n$. Recall the relationship between the monomial generators and the Young idempotent generators for $S(S^1 \times D^2)$, we can rewrite the monomial basis elements of the subspace $C_{-m_{1}} \times C_{n_1}$ in terms of $Q_{\nu, \sigma}$ with $|\nu| = |m_{1}|$ and $|\sigma| = |n_1|$. The result follows.
\end{pf}

\begin{Cor}
When $|\lambda| \geq 1$, $|\mu| \geq 1$, $(1-x^{2(|\mu|-|\lambda|)}v^{2}s^{-2cn(c)})Q_{\lambda, \mu}$ can be written as a linear combination of the $Q_{\nu, \sigma}$'s with $(\nu, \sigma) \in \varGamma(\lambda, \mu)$, where $c$ is any extreme cell of $\lambda$. Thus if $\lambda$ is not a rectangular Young diagram, $Q_{\lambda, \mu}$ can be written as a linear combination of the $Q_{\nu, \sigma}$'s with $(\nu, \sigma) \in \varGamma(\lambda, \mu)$.
\end{Cor}
\begin{pf} 
Simplify Equation (5.3), 
$$
Q_{\lambda, \mu}=\raisebox{-25mm}{\ \epsfxsize1.8in\epsffile{ylamu1re1a.ai}}\ \equiv \ \raisebox{-25mm}{\ \epsfxsize1.8in\epsffile{ylamu1re2a.ai}}
$$
By applying the skein relations in a similar computation to that in the proof of Lemma 5.3,
here we use Corollary 3 in place of Proposition 3.4. 
$$\raisebox{-25mm}{\ \epsfxsize1.8in\epsffile{ylamu1re2a.ai}}=x^{2(|\mu|-|\lambda|)}v^{2}s^{-2cn(c)}\raisebox{-25mm}{\ \epsfxsize1.8in\epsffile{ylamu1re1a.ai}}
$$
$$
+\sum_{(\nu, \sigma) \in \varGamma(\lambda, \mu)}b_{\nu\sigma}Q_{ \nu, \sigma}.
$$
By a similar proof to Lemmas 5.3, the elements of the form 
$$\raisebox{-20mm}{\ \epsfxsize1.5in
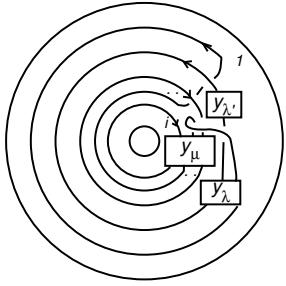
\epsffile{ymulai.ai}}$$
can be written as a linear combination of the $Q_{\nu, \sigma}$ with $(\nu, \sigma) \in \varGamma(\lambda, \mu)$. Substituting these into Equation (5.3), we have:
 $$(1-x^{2(|\mu|-|\lambda|)}v^{2}s^{-2cn(c)})Q({\lambda, \mu})\equiv \sum_{(\nu, \sigma) \in \varGamma(\lambda, \mu)}b_{\nu\sigma}Q_{ \nu, \sigma}$$.
\end{pf}

Combining Lemma 5.2 and Corollary 6, we have:
\begin{Cor}
$(v^{-2}s^{2cn(c)}-v^{2}s^{-2cn(c')})Q({\lambda, \mu})\equiv \sum_{(\nu, \sigma) \in \varGamma(\lambda, \mu)}d_{\nu\sigma}Q_{ \nu, \sigma}$.
\end{Cor}
\noindent Here $d_{\nu\sigma} \in k$, $c$ and $c'$ are extreme cells of $\mu$ and $\lambda$, respectively.

\begin{Lem}
Every generator $Q_{\lambda, \mu}$ with $|\lambda|\ne |\mu|$ in $S(S^1 \times D^2)$ is torsion in $S(S^1 \times S^2)$. i.e. $S$ is a torsion submodule.
\end{Lem}
\begin{pf}
We will proceed by induction on $|\lambda|$.

Lemma 5.1 proves the result for $|\lambda|=0$. We now consider the case $|\lambda|\geq 1$. When $|\nu|\ne |\sigma|$, suppose $Q_{\nu, \sigma}$ is torsion for $|\nu|<|\lambda|$. By Lemma 5.2 (or Corollary 7) and the induction hypothesis, the result follows.
\end{pf}

\begin{Lem}
For $Q_{\lambda, \mu}$ with $|\lambda|= |\mu|$, we have $t_1Q_{\lambda, \mu}\equiv t_2\phi$ in $S(S^1\times S^2)$, where $t_1, t_2$ are two scalars in $k$ and $\phi$ is the empty link.
\end{Lem}
\begin{pf}
We will proceed by induction on $|\lambda|$ as well.
(1) $Q_{0,0}=\phi$ is the trivial case. 

(2) For $|\lambda|>0$, suppose the result holds for $Q_{\nu, \sigma}$ for all $|\nu|<|\lambda|$. By the remark following Lemma 5.2, we have $$(1-v^{-2}s^{2cn(c)})Q_{\lambda,\mu}\equiv \sum_{(\nu, \sigma) \in \varGamma(\lambda, \mu)}a_{\nu\sigma}Q_{\nu, \sigma}.$$
where $a_{\nu\sigma} \in k$. By the induction hypothesis, the result follows.
\end{pf}

\subsection{Proof of Theorem 2 \& Proposition 1.3}

\begin{pf}
Recall that $S(S^1\times D^2)=S_0 \oplus S$, where $S_0$ is generated by \{$Q_{\lambda, \mu}\mid |\lambda|=| \mu|$\}, and $S$ is generated by \{$Q_{\lambda, \mu}\mid |\lambda|\ne| \mu|$\}.
(i) In the quotient space $S(S^1\times S^2)/<\phi>$, note $<\phi>=0$. So by Lemma 5.5, all generators of $S_0$ are torsion. From Lemma 5.4, $S$ is a torsion submodule in $S(S^1\times S^2)$. Therefore $S(S^1\times S^2)/<\phi>$ is torsion.

(ii) If all elements of the form $s^{2n}-v^{2}$ are invertible in $k$ for $n \in {\Bbb Z}$, then each $Q_{\lambda, \mu}\in S_{0}$ is a scalar multiple of $\phi$. i.e. $S_0=<\phi>$. This proves Proposition 1.3. Again $S$ is a torsion submodule from Lemma 5.4. The result $S(S^{1} \times S^2)={\hbox {k-torsion module}}\ \oplus <\phi>$ follows.

(iii) In addition, when we make the assumption that all elements of the form $s^{2n}-v^{4}$ are invertible in $k$ for $n \in {\Bbb Z}$, then $S_0=<\phi>$ from (ii); $S={0}$ which follows by induction and Corollary 7. Therefore $S(S^{1} \times S^2)=<\phi>$.
\end{pf}

%\end{document}
\section{Proof of Theorem 3}
\begin{Lem}
Let $f(x,v,s)=\sum_{(i,j,k)}a_{ijk}x^{i}v^{j}s^{k}$ be a polynomial in $x, v, s$ with nonzero coefficients $a_{ijk}$. If $n$ is sufficiently large, then all the exponents $i+jn^2-kn$ of $x$ in $f(x,x^{n^2},x^{-n})=\sum_{(i,j,k)}a_{ijk}x^{i+jn^2-kn}$ are distinct.
\end{Lem}
\begin{pf}
Now take any two distinct pairs $(i_1,j_1,k_1)$ and $(i_2,j_2,k_2)$ from the finite collection $\{(i,j,k)\}$ of the exponent of $x^{i+jn^2-kn}$ in $f(x,x^{n^2},x^{-n})$. Let $Y=(i_{1}+j_{1}n^2-k_{1}n)-(i_{2}+j_{2}n^2-k_{2}n)$. So $Y=(j_{1}-j_{2})n^2 -(k_{1}-k_{2})n + (i_{1}-i_{2})=An^2+Bn+C$, where we denote the corresponding coefficients by $A, B, C$. Note $A, B, C$ are bounded and are not all zero. 

(1) If $A\ne 0$, we have a finite collection of quadratic functions in the variable $n$ when the pairs vary over all $(i,j,k)$. These quadratic functions only have finitely many bounded roots. So we can take integer $D$ such that if $n>D$, $Y=An^2+Bn+C\ne 0$. i.e. $i_{1}+j_{1}n^2-k_{1}n \ne i_{2}+j_{2}n^2-k_{2}n$.

(2) If $A=0$, i.e. $j_{1}=j_{2}$, then $B, C$ are not both zero. We can take integer $E$ such that if $n>E$, $Bn+C\ne 0$ for the finite collection of $B, C$. i.e. $i_{1}+j_{1}n^2-k_{1}n \ne i_{2}+j_{2}n^2-k_{2}n$ when $j_{1}=j_{2}$.

Then the exponents $i+jn^2-kn$ are distinct if $n>\max\{D, E\}$.
\end{pf}

{\bf Proof of Theorem 3}:
If $k$ is a subring of the field of rational functions in $x, v, s$ which contains $x^{\pm 1}, v^{\pm 1}, s^{\pm 1}$ over ${\Bbb {C}}$. We want to show the submodule $<\phi>$ is free in $S(S^1\times S^2)$ over $k$. 

Let $\varLambda={\Bbb C}[x^{\pm 1}, v^{\pm 1}, s^{\pm 1}]$. Note $\varLambda \subseteq k$. Let $\varLambda_{0}={\Bbb C}[x, v, s]$. Let $S_{\varLambda}(M)$ denote the Homflypt skein module of $M$ over $\varLambda$. 

By Corollary 2 in section 2, $S(S^1\times S^2)\cong  S(S^1\times D^2)/ R$. Suppose that $<\phi>$ is not a free submodule of $S(S^1\times S^2)$, then there exists a nonzero rational function $\dfrac{f(x,v,s)}{g(x,v,s)}$ such that $\dfrac{f(x,v,s)}{g(x,v,s)}\phi=\sum_{i}\dfrac{f_{i}(x,v,s)}{g_{i}(x,v,s)}R_{i} $, where the sum is a finite sum, $f, f_{i},g, g_{i} \in \varLambda_{0}$ and $R_{i}\in R$. We can clear the denominators and have $P(x,v,s)\phi=\sum_{i}P_{i}(x,v,s)R_{i}$, where $P, P_{i}\in \varLambda_{0}$. Note $P(x,v,s)$ is nonzero. Therefore $P(x,v,s)\phi=0$ in $S_{\varLambda_{0}}(S^1\times S^2)$. 

Let $M$ be an oriented $3$-manifold. For each pair of integers $N, K$ with $N\geq 2$, $K \geq 1$, choose a primitive $2N(N+K)$th root of unity denoted by $t$. $t$ depends on $N, K$. Y. Yokota  \cite{Y97} has shown that there is a ${\Bbb C}$-valued invariant:
$$I_{N,K}(M)=\theta<\theta\Omega_{K}>_{U_{-}}^{\sigma}<\theta_{\Omega_{K}}, \cdots, \theta_{\Omega_{K}}>_{D}$$
where $\theta \in {\Bbb C}$, $\theta^{-2}=<\Omega_{K}>_{U}$, and $D$ is a diagram in $S^2$ which is a surgery description for $M$ and $\sigma$ is the signature of the linking matrix of $D$. Also see \cite{L97}. Note $I_{N,K}(S^1 \times S^2)=1$. Let $L$ be a framed oriented link in $M$. In the usual way, one may extend $I_{N,K}(M)$ to an invariant:
$$I_{N,K}(M, L)=\theta<\theta\Omega_{K}>_{U_{-}}^{\sigma}<\theta_{\Omega_{K}}, \cdots, \theta_{\Omega_{K}, \alpha, \cdots, \alpha}>_{D\cup L}$$
here $\alpha$ represents the core of the annulus with the trivial framing, and $\alpha$ occurs in the positions corresponding to the components of $L$. Note $I_{N,K}(S^1 \times S^2, \phi)=I_{N,K}(S^1 \times S^2)=1$. 

We now give ${\Bbb C}$ a $\varLambda$-module structure where: $x$ acts on ${\Bbb C}$ as $t^{-1}$, $v$ acts on ${\Bbb C}$ as $t^{-N^2}$, $s$ acts on ${\Bbb C}$ as $t^{N}$. We denote ${\Bbb C}$ with a $\varLambda$-module structure by ${\Bbb C}_{N,K}$. Let ${\cal I}_{N,K}(S^1 \times S^2)$ be the map: $\varLambda_{0}F(S^1\times S^2) \to {{\Bbb C}_{N,K}}$ given by ${\cal I}_{N,K}(S^1 \times S^2)(\sum_{i}f_{i}(x,v,s)L_{i})=\sum_{i}f_{i}(t^{-1},t^{-N^2},t^{N})I_{N,K}(S^1 \times S^2)(L_{i})$.
Note the submodule $\varLambda_{0}H(S^1\times S^2)$ of $\varLambda_{0}F(S^1\times S^2)$ is mapped to zero by ${\cal I}_{N,K}(S^1 \times S^2)$ due to the skein relation of Yokota. So we have an induced homomorphism:
$$Z_{N, K}:\ S_{\varLambda_{0}}(S^1\times S^2) \to {\Bbb C}_{N,K}.$$
Note the evaluation $Z_{N, K}(\phi)=1$. Therefore $P(t^{-1},t^{-N^2},t^{N})=0$ for all positive integers $K \geq 1$ and $N\geq 2$. 

Assume $P(x,v,s)=\sum_{(i,j,k)}a_{ijk}x^{i}v^{j}s^{k}$, where the coefficients $a_{ijk}$ are nonzero. By the previous lemma, we can take an $n$ such that $P({\tau}^{-1},{\tau}^{-n^2},{\tau}^{n})$

$=\sum_{(i,j,k)}a_{ijk}{\tau}^{-(i+jn^2-n)}$ is nonzero, but by the above $P({\tau}^{-1},{\tau}^{-n^2},{\tau}^{n})$ has the infinite collection of roots $t$ for $K \geq 1$ and the chosen $n$. This is a contradiction.%\end{document}

%\begin{thebibliography}{10}
%\end{thebibliography}

\end{document}